\documentclass[reqno]{amsart}

\usepackage{tikz-cd}
\usepackage{thmtools}

\makeatletter
\@ifundefined{newcounteralias}{}{%
  \renewcommand\thmt@autorefsetup{%
    \@xa\def\csname\thmt@envname autorefname\@xa\endcsname
      \@xa{\thmt@thmname}%
  }%
}
\makeatother

\usepackage{enumitem}
\usepackage{quiver}
\usepackage[most]{tcolorbox}

\usepackage{csquotes} 
\usepackage[style=alphabetic, backend=biber, sorting=nyt, doi=false,
isbn=false, url=false, hyperref, maxnames=20,
maxalphanames=3]{biblatex}

\newbibmacro{string+doiurl}[1]{%
  \iffieldundef{doi}{%
    \iffieldundef{url}{%
          #1%
    }{%
      \href{\thefield{url}}{#1}%
    }%
  }{%
    \href{https://dx.doi.org/\thefield{doi}}{#1}%
  }%
}
\DeclareFieldFormat{title}{\usebibmacro{string+doiurl}{\mkbibemph{#1}}}
\DeclareFieldFormat[article,incollection,inproceedings,unpublished,misc,book]{title}%
    {\usebibmacro{string+doiurl}{\mkbibquote{#1}}}
\DeclareFieldFormat[online]{title}{\mkbibquote{#1}}

\DeclareFieldFormat{titlecase}{\MakeTitleCase{#1}}
\newrobustcmd{\MakeTitleCase}[1]{%
  \ifthenelse{\ifcurrentfield{booktitle}%
    \OR\ifcurrentfield{journaltitle}%
    \OR\ifentrytype{book}%
    \OR\ifentrytype{booklet}
    \OR\ifentrytype{collection}%
    \OR\ifentrytype{proceedings}%
         \AND\NOT\iffieldequalstr{langid}{french}}
    {#1}
    {\MakeSentenceCase{#1}}}

\renewbibmacro{in:}{%
  \ifentrytype{article}{}{\printtext{\bibstring{in}\intitlepunct}}}

\AtBeginBibliography{\small}

\usepackage[]{mdframed}
\usepackage{upgreek}
\usepackage{mathrsfs}

\usepackage{numprint}

\usepackage{hyperref}
\usepackage[capitalize]{cleveref}

\DeclareFieldFormat{postnote}{#1}
\DeclareFieldFormat{multipostnote}{#1}

\usepackage{enumitem}
\setlist[enumerate]{label=\arabic*.}

\crefformat{equation}{(#2#1#3)}
\Crefformat{equation}{(#2#1#3)}

\newtheorem{theorem}{Theorem}[section]
\newtheorem{corollary}[theorem]{Corollary}
\newtheorem{lemma}[theorem]{Lemma}
\newtheorem{proposition}[theorem]{Proposition}

\theoremstyle{definition}
\newtheorem{remark}[theorem]{Remark}
\newtheorem{definition}[theorem]{Definition}
\newtheorem{example}[theorem]{Example}

\DeclareMathOperator{\ord}{ord}

\newcommand{\Z}{\mathbb{Z}}

\newcommand{\C}{\mathbb{C}}
\newcommand{\Q}{\mathbb{Q}}
\newcommand{\Proj}{\mathbb{P}}
\newcommand{\UH}{\mathbb{H}}

\newcommand{\Pa}{\mathcal{P}}
\newcommand{\F}{\mathcal{F}}
\newcommand{\OO}{\mathcal{O}}

\newcommand{\Sg}{\mathscr{S}}

\newcommand{\f}{\mathfrak{f}}
\newcommand{\m}{\mathfrak{m}}

\newcommand{\SL}{\mathrm{SL}}
\newcommand{\GL}{\mathrm{GL}}
\newcommand{\PSL}{\mathrm{PSL}}
\newcommand{\Sch}{\mathrm{Sch}}
\newcommand{\Sets}{\mathrm{Sets}}
\DeclareMathOperator{\Spec}{Spec}
\DeclareMathOperator{\Gal}{Gal}
\DeclareMathOperator{\Hom}{Hom}
\DeclareMathOperator{\Sym}{Sym}

\newcommand{\crt}{\mathrm{crt}}

\newcommand{\domain}{\mathrm{domain}}
\newcommand{\codomain}{\mathrm{codomain}}

\newcommand{\Stab}{\mathrm{Stab}}
\renewcommand{\mod}{\ensuremath{\ \text{mod}\ }}

\newenvironment{numberedcenter}
  {\begin{equation}\begin{minipage}{0.8\textwidth}\centering}
  {\end{minipage}\end{equation}}

\title[Modular polynomials for genus 0 covers of $X(1)$]{Modular polynomials for genus 0 covers of $X(1)$ have smaller height}
\author{Ilan Ehrlich and Jean Kieffer}

\begin{document}

\begin{abstract}
The classical modular polynomials classify pairs of $\ell$-isogenous elliptic curves via their $j$-invariant. Following work of Atkin, Morain, Enge and others, it has been observed experimentally that using invariants corresponding to richer level structures, such as the Weber modular function $\f$, gives rise to modular polynomials of significantly smaller height. We provide a proof of this general phenomenon: given an invariant on a genus 0 modular curve which is a degree-$d$ cover of $X(1)$, we show that the height of the associated $\ell$-modular polynomial is bounded above by $\frac{6}{d} \ell \log \ell + O(\ell \sqrt{\log \ell})$. Along the way, we recast the construction of such modular polynomials within the modern framework of modular curves.
\end{abstract}

\maketitle
\section{Introduction}

Let~$\ell$ be a prime number. The well-known classical modular polynomial $\Phi_\ell(x,y)$ is a bivariate polynomial with integer coefficients, of degree $\ell+1$ and monic in both variables, and satisfying the following property: if~$k$ is a field of characteristic different from $\ell$, and $E, E'$ are elliptic curves over~$k$, then $E$ and~$E'$ are $\ell$-isogenous over an extension of~$k$ if and only if $\Phi_\ell(j(E),j(E')) = 0$. As a universal algebraic formula encoding the presence of $\ell$-isogenies, the classical modular polynomial has numerous algorithmic applications, for instance \cite{schoofCountingPointsElliptic1995,sutherlandComputingHilbertClass2011,bissonComputingEndomorphismRing2011}.

Ultimately, the cost of evaluating~$\Phi_\ell$ and of using it in algorithms is proportional to the size of its coefficients \cite{engeComputingModularPolynomials2009,brokerModularPolynomialsIsogeny2012,sutherlandEvaluationModularPolynomials2013}. Let~$h(\Phi_\ell)$ denote the height of~$\Phi_\ell$, defined as $\log\max|c|$ as $c$ runs through the coefficients of~$\Phi_\ell$. Cohen~\cite{cohenCoefficientsTransformationPolynomials1984} proved that
\[
h(\Phi_\ell) \sim_{\ell\to+\infty} 6\ell\log \ell.
\]
See also \cite{brokerExplicitHeightBound2010,breuerExplicitBoundsCoefficients2025} for explicit upper bounds. As a concrete example, the largest coefficient of~$\Phi_\ell$ for $\ell=101$ has 1731 decimal digits \cite{sutherlandDatabaseModularPolynomials}.

It has soon been observed experimentally that modular polynomials for other modular functions than~$j$ can be significantly more compact, leading to improved performances to compute isogenies between elliptic curves in practice: see for instance \cite[§5.2]{Bröker_2011}. Perhaps the most well-known example of such an alternative invariant is the Weber modular function~$\f$, which satisfies the following degree~72 relation to the $j$-invariant:
\[
(\f^{24}-16)^3-j\f^{24}=0.
\]
The attached modular polynomial $\Phi_\ell^\f$ is then expected to satisfy
\begin{equation}
h(\Phi_\ell^\f) \underset{\ell\to+\infty}{\sim} \frac{1}{72}h(\Phi_\ell).
\label{eq:height-weber}
\end{equation}
In addition, $\Phi_\ell^\f$ is sparse, as only one coefficient in 24 can be nonzero. Consequently, $\Phi_\ell^\f$ should be 1728 times more compact to write than~$\Phi_\ell$ as a bivariate polynomial \cite{Elkies1997EllipticAM,sutherland_weber_modular_polynomials} and \cite[§7.3]{brokerModularPolynomialsIsogeny2012}. The practical gain is evident.

Recently, A.~Zhang \cite{zhang2026explicitheightboundsmodular} obtained an upper bound on~$h(\Phi_\ell^\f)$ 
based on the explicit $q$-expansion of $\f$, from which~\eqref{eq:height-weber} can be deduced. Our primary goal in this paper is to provide an independent upper bound that not only applies to the Weber function, but more generally to any invariant arising from a genus zero modular cover of the $j$-line. This result ultimately provides a better understanding of the asymptotic complexity of any algorithm involving these polynomials.

\subsection*{Main results} 

Recall that the group $\SL_2(\mathbb{Z})$ acts on the extended upper half plane~$\UH^* = \UH \cup \mathbb{P}^1(\Q)$ by Möbius transformations. The quotient $X(1) = \SL_2(\Z)\backslash \UH^*$ is a smooth, projective algebraic curve over~$\mathbb{C}$, isomorphic to~$\smash{\mathbb{P}^1_\Q}$ via the $j$-invariant. If~$k$ is an algebraically closed field, such as~$\overline{\Q}$, then $X(1)(k)$ is naturally identified with the set of isomorphism classes of elliptic curves over~$k$.

If~$F$ is a polynomial with coefficients in~$\overline \Q$, its (absolute logarithmic) height~$h(F)$ is defined as in \cite[§B.7]{hindry2000diophantine}. This notion of height specializes to the previous one in the case of polynomials with integer coefficients.

\begin{theorem}\label{thm:main-height}
    Let $X$ be a smooth, projective curve of genus zero defined over a number field~$K$ endowed with a covering map $p: X\to X(1)$ of degree~$d$, and choose an isomorphism 
    $f:X \overset{\sim}{\longrightarrow} \mathbb{P}^1_{K}$. Let~$\ell$ be a prime number. Suppose that there exists a polynomial $\Phi_\ell^f \in \overline \Q[x,y]$ such that:
    \begin{enumerate}
        \item \label{intro degree condition}
        $\Phi_\ell^f$ has degree $\ell+1$ in both variables,
        and is monic as a polynomial in $y$.

        \item \label{intro isogenous condition}
        For every pair $(Q,Q')\in 
        X(\overline \Q) \times X(\overline \Q)$ 
        satisfying
        \[
            \Phi_\ell^f(f(Q),f(Q')) = 0,
        \]
        the pair $(p(Q),p(Q'))$ is either $(\infty,\infty)$ or corresponds to~$\ell$-isogenous elliptic curves.
    \end{enumerate}
    Then we have
    \[
    h(\Phi_\ell^f) \leq \frac{6}{d}\ell\log\ell + O_{p,f}\left(\ell\sqrt{\log(\ell)}\right).
    \]
    The implied constants in this estimate are independent of~$\ell$ and depend only on~$p,f$.
\end{theorem}

The conditions~\eqref{intro degree condition}--\eqref{intro isogenous condition} of \cref{thm:main-height} are arguably minimal requirements for a well-behaved modular polynomial of level~$\ell$ relative to the invariant~$f$. Our proof is based on the relation between the height of a polynomial and the height of its values at well-chosen points (in our case, roots of unity), and only relies of those few properties as opposed to a specific construction. See also \cref{cor:explicit}, which provides explicit upper bounds in terms of the map $j\circ p\circ f^{-1}: \mathbb{P}^1\to \mathbb{P}^1$.

It is not a priori obvious why polynomials satisfying the conditions of \cref{thm:main-height} should exist for a given invariant~$f$ and prime level~$\ell$, even if~$X$ is itself a connected modular curve, i.e.~a quotient $\Gamma\backslash \UH^*$ for some congruence subgroup $\Gamma\subseteq \SL_2(\Z)$. For example, if we take~$f=\f$ to be the Weber function, a direct transcription of the existence proof for the classical modular polynomials~$\Phi_\ell$ in \cite[§11]{Cox1997} would fail, as the corresponding~$X$ has two cusps. To the best of our knowledge, a similar proof of the existence of modular polynomials for the Weber function has never appeared in print. Our second main goal is to fill this gap and explain why modular polynomials exist for a large class of invariants and levels~$\Gamma$, complementing partial results from
\cite[§5.1]{Bröker_2011}.

To answer this constructive question, we use further facts on modular curves (see also the recent preprint~\cite{onukiGeometricConstructionModular2026} for a different route).
In the language of §\ref{section: modular curves}, let~$n\geq 1$ be an integer, let~$H$ be a subgroup of~$\GL_2(\Z/n\Z)$, and let~$X_H$ denote the associated modular curve, seen as an algebraic curve over~$\Q$. It is a coarse moduli space for elliptic curves equipped with a level~$H$ structure. 
If~$\pi_n: M_2(\Z)\to M_2(\Z/n\Z)$ denotes the natural projection, we define
\begin{equation}
    \label{eq:GammaH}
    \Gamma_H = \pi_n^{-1}(H)\cap \SL_2(\Z).
\end{equation}
Every congruence subgroup of~$\SL_2(\Z)$ arises in this way. The curve~$X_{\Gamma_H} = \Gamma_H\backslash \UH^*$ is then identified with a geometric connected component of~$X_H$, and is defined over the number field 
\begin{equation}\label{eq:def KH}
    K_H = \Q(\zeta_n)^{\det H}.
\end{equation}
For any~$\tau\in \UH$, the complex elliptic curve $E_\tau = \C/(\tau\Z+\Z)$ carries a canonical level~$H$ structure that only depends on the class of~$\tau$ modulo~$\Gamma_H$. For this to hold, we add a transposition to the traditional definition of a level~$H$ structure.

The geometric connected components of~$X_H$ are indexed by elements of the group $(\Z/n\Z)^\times/\det H$, depending on the value of the Weil pairing on the basis elements of the level~$H$ structure, and taking $\ell$-isogenies has the effect of taking the $\ell$-th power of this Weil pairing. Therefore, a natural condition to obtain modular polynomials on the single component $X_{\Gamma_H}$, as in \cref{thm:main-height}, is to assume that $(\ell \text{ mod } n)\in \det H$. 

For each prime~$\ell$ such that $(\ell\text{ mod }n)\in \det H$, we construct in §\ref{subsec:complex-isog} an integral matrix~$ M_\ell$, of determinant~$\ell$, such that the double coset $\Gamma_H\backslash (\Gamma_H M_\ell \Gamma_H)$ contains exactly $\ell+1$ elements. Then, we prove:

\begin{theorem}
\label{thm:main-existence-1} Let~$n\geq 1$, let~$H$ be a subgroup of $\GL_2(\Z/n\Z)$, and define~$\Gamma_H$ and~$K_H$ as above. Assume that~$X_{\Gamma_H}$ has genus zero, and choose an isomorphism $f: X_{\Gamma_H} \to \mathbb{P}^1_K$ such that $f(\infty)=\infty$. Then for any prime number $\ell\nmid n$ such that $(\ell\text{ mod } n)\in \det H$, there exists a rational fraction $\Phi_\ell^f(x,y)\in K_H(x)[y]$ with the following properties.
\begin{enumerate}
    \item The denominator of $\Phi_\ell^f$ may only vanish at the images of cusps of~$X_{\Gamma_H}$ under~$f$.
    \item For every~$\tau\in \UH$, we have
\[
\Phi_\ell^f(f(\tau),y) = \prod_{\gamma\in \Gamma_H\backslash (\Gamma_H M_\ell \Gamma_H)} (y - f(\gamma\tau)).
\]
    \item For any $\tau,\tau'\in \UH$, we have $\Phi_\ell^f(f(\tau),f(\tau')) = 0$
if and only if there exists an $\ell$-isogeny $\phi: E_\tau \to E_\tau'$ that respects the canonical level~$H$ structures.
    \item The degree of $\Phi_\ell^f$ in $y$ and the degree of the numerator of $\Phi_\ell^f$ in $x$ are both equal to $\ell+1$. In addition, $\Phi_\ell^f$ is monic in $y$.
\end{enumerate}
\end{theorem}
The rational fraction~$\Phi_\ell^f$ will be denominator-free precisely when its coefficients have no poles as rational fractions of~$x$. This property only depends on the action of the double coset $\Gamma_H\backslash (\Gamma_H M_\ell \Gamma_H)$ on the cusp $\infty$:

\begin{theorem}
\label{thm:main-existence-2}
    Under the assumptions of \cref{thm:main-existence-1}, the rational fraction $\Phi_\ell^f$  is denominator-free if and only if for each coset $\gamma\in \Gamma_H\backslash (\Gamma_H M_\ell \Gamma_H)$ and each cusp $\kappa\in \mathbb{P}^1(\Q)\backslash (\Gamma_H \infty)$, we have $\gamma\kappa\notin \Gamma_H \infty$.
\end{theorem}

We further show that the condition on cusps in \cref{thm:main-existence-2} is equivalent to an explicit group-theoretic condition involving only $H$ and $(\ell \mod n)$. When~$\Phi_\ell^f$ is a bivariate polynomial, one may ask whether it is symmetric in~$x$ and~$y$, as the classical modular polynomial~$\Phi_\ell$; we are also able to answer that question.

\begin{theorem}\label{thm: summary conditions}
    Keep the assumptions of \cref{thm:main-existence-1}, and let~$u = (\ell \text{ mod } n)$. Then:
    \begin{enumerate}
        \item \label{item: intro denominator}The condition on cusps in \cref{thm:main-existence-2} holds if and only if~$H$ contains matrices of the form
        \[\varepsilon_1 \begin{pmatrix}
            u^{-1} & v \\ 0 & u
        \end{pmatrix} 
        \quad\text{and}\quad
        \varepsilon_2 \begin{pmatrix}
            1 & w \\ 0 & u
        \end{pmatrix}
        \]
        for some $\varepsilon_1,\varepsilon_2\in \{\pm 1\}$ and $v,w\in \Z/n\Z$.
        \item \label{item: intro symmetry} We have $\Phi_\ell^f(x,y) = \Phi_\ell^f(y,x)$ if and only if, in addition, $\varepsilon uI\in H$ for some $\varepsilon\in \{\pm 1\}$.
    \end{enumerate}
\end{theorem}

Finally, under the assumptions of \cref{thm: summary conditions}, we present sufficient conditions for the polynomial $\Phi_\ell^f$ to have integer coefficients: see \cref{cor: coefficients in Z}.

As an application, we verify that the conditions of \cref{thm: summary conditions,cor: coefficients in Z} hold in two folklore cases, namely that of the modular function $\gamma_2$ (the cube root of~$j$) and the Weber function $\f$. The associated modular polynomials are thus symmetric polynomials in $\Z[x,y]$, and \Cref{thm:main-height} provides an upper bound on their heights.

Note that the Weber invariant $\f$ is optimal in the sense that its associated modular curve is the genus~$0$ modular cover of~$X(1)$ of largest possible degree, according to the classification from \cite{CHUA2004408}. However, \cref{thm:main-height} does not require~$X$ to be a connected modular curve, and one might ask whether there exist other covers of~$X(1)$ of degree higher than 72, together with polynomials satisfying its hypotheses. A positive answer would produce polynomials classifying isogenous elliptic curves with smaller heights than those known to this day.
 
\subsection*{Organization of the paper}
The paper is organized as follows. In~§\ref{section: modular curves}, we provide general background on modular curves; see also Appendix~\ref{appendix} for a detailed presentation of modular curves as moduli spaces. We use this framework in §\ref{sec:modular curves and isogenies} to produce an analytic uniformization of $\ell$-isogenies respecting level structures, which arises from morphisms between modular curves. This construction is crucial to our study of modular polynomials in §\ref{section: construction}, where we prove \cref{thm:main-existence-1,thm:main-existence-2,thm: summary conditions}.
Section §\ref{section heights} is devoted to the proof of \cref{thm:main-height}. Finally, we apply our results to the examples of  the modular function $\gamma_2 = j^{1/3}$ and the Weber modular function $\f$ in §\ref{section: examples}.

\subsection*{Notation}
Throughout this article, $n$ denotes an arbitrary positive integer,
and $\ell$ is a prime number not dividing $n$.
We denote by $\pi_n$ the reduction map modulo~$n$, whether applied
to an integer, a class modulo a multiple of $n$, or a matrix with
integer or modular entries. We also set
$u:=\pi_n(\ell)\in(\Z/n\Z)^\times$.

Concerning modular curves, our notation is as follows. We denote by $H$ the subgroups of $\GL_2(\Z/n\Z)$ under consideration.
To such a subgroup $H$ we associate the subgroup $\Gamma_H\leq \SL_2(\Z)$ as in \cref{eq:GammaH} and the number field $K_H$ as in \cref{eq:def KH}. The modular curve associated to $H$ is denoted by $Y_H$ and its compactification by $X_H$; these are algebraic curves over~$\Q$. Their geometric connected components are respectively denoted by $Y_H^{(\beta)}/K_H$ and $X_H^{(\beta)}/K_H$, where $\beta\in (\Z/n\Z)^\times/ \det H$. The Riemann surface $X_{\Gamma_H}= \Gamma_H\backslash \UH^*$ is then identified as a complex curve to $X_H^{(-1)}$. 
In §\ref{section: construction}, we make the assumption on $H$
that $X_{H}^{(-1)}$ has genus~$0$, and denote by
$f:X_{H}^{(-1)}\to\Proj_{K_H}^1$ an isomorphism such that
$f(\infty)=\infty$ over $\C$. 

We also use the following notation for matrices.
The expression $D(a,b)$ stands for the diagonal matrix with entries
$a,b$, whether these are rational numbers or modular classes.
As usual, we write
\[
S=\begin{pmatrix}
    0 & -1 \\ 1 & 0
\end{pmatrix}
\qquad\text{and}\qquad
T=\begin{pmatrix}
    1 & 1 \\ 0 & 1
\end{pmatrix}
\]
for the standard matrices in $\SL_2(\Z)$ or $\SL_2(\Z/n\Z)$.
Starting from §\ref{subsec:isog-components}, we add the assumption that
$u\in\det H$. The matrices~$M_u$, $\gamma_u$, $M_\ell$ are defined  after \cref{lem:conditions diag}.
Finally, \cref{equivalence} introduces the matrices
\[
U(u,v):=\begin{pmatrix}
    u^{-1} & v \\ 0 & u
\end{pmatrix},
\qquad
M(u,w):=\begin{pmatrix}
    1 & w \\ 0 & u
\end{pmatrix}
\]
in $\GL_2(\Z/n\Z)$.
\subsection*{AI disclosure} ChatGPT 6 Sol was used to prove \cref{prop:galois-sheaf}. Claude Opus 5.5 was used in preliminary discussions about \cref{prop:moduli-product,thm:algebraic diag}.

\subsection*{Acknowledgements} The first author was partially supported by the ANR grant ANR-21-CE39-0008-01. He wishes to express his warm gratitude to Inria and the CARAMBA team, as well as to the University of Luxembourg. He is also particularly grateful to Antonella Perucca.

\section{Background on modular curves}\label{section: modular curves}

In this section, we provide general background on modular curves. This material will then be used in §\ref{sec:modular curves and isogenies} in relation with $\ell$-isogenies. 

For simplicity, we work over fields in this presentation, and refer to Appendix~\ref{appendix} for the full characterization of modular curves as moduli spaces of elliptic curves over schemes. This full characterization will be logically necessary to assert that certain maps between modular curves in §\ref{sec:modular curves and isogenies} are in fact algebraic morphisms, but is perhaps less accessible. We also choose to work in characteristic zero for simplicity.

\subsection{Level structures over fields}\label{subsec: level struct for fields}
Let~$n\geq 1$, and let $E$ be an elliptic curve over a field $k$ of characteristic zero. Recall that $E[n](\overline k) \simeq (\Z/n\Z)^2$ as an abstract group, and $E[n]$ is fully described by the data of $E[n](\overline k)$ together with the Galois representation
    $$
    \rho_{E,n}:\Gal(\overline k/k) \to \mathrm{Aut}(E[n](\overline k)) \simeq \GL_2(\Z/n\Z),
    $$
this last isomorphism being only well-defined up to conjugacy.

\begin{definition}
\label{def:level-structure}
A \textit{full level~$n$ structure} on~$E$ over~$k$ is an isomorphism of groups
$$E[n](k) \to (\Z/n\Z)^2.$$
Such a structure exists if and only if each point of~$E[n]$ is defined over~$k$, i.e.~$\rho_{E,n}$ is trivial.

The group $\GL_2(\Z/n\Z)$ acts (on the right) on the set of full level structures~$\alpha$ on~$E$ over~$\overline k$ via
\begin{equation}\label{eq:right action}
(\alpha, h) \mapsto h^t \circ \alpha.
\end{equation}
Fixing a subgroup $H\leq \GL_2(\Z/n\Z)$, a \textit{level~$H$ structure} on $E$ over~$k$ is a class of full level structures $\alpha$ on~$E$ over~$\overline k$ modulo the action of~$H$ as above, that is compatible with the Galois representation: for every $\sigma\in \Gal(\overline k/k)$, there must exist $h_\sigma\in H$ such that 
$$\alpha \circ \rho_{E,n}(\sigma) = h_\sigma^t \circ \alpha.$$
Such a structure exists if and only if the image of $\rho_{E,n}$ is contained in~$H$ up to conjugacy.
\end{definition}

\begin{remark}
\label{rem:lmfdb}
    This definition of level~$H$ structures is not in accordance with the traditional one, used in the LMFDB \cite{thelmfdbcollaborationLfunctionsModularForms2026} for instance, which uses the left action $(g,\alpha)\mapsto g\circ \alpha$ instead of~\eqref{eq:right action}. We chose the above convention as it interacts better with isomorphisms between complex elliptic curves in~§\ref{subsec:complex} below. The modular curve~$X_H$ used in this paper will therefore correspond to level $H^t$ structures in the traditional sense; note that $H^t$ is not conjugate to~$H$ in general.
\end{remark}

From now on, we will identify level~$H$ structures over fields with one of their representative isomorphisms from $E[n](\overline k)$ to $(\Z/n\Z)^2$; this representative is well-defined up to the action of~$H$. Two pairs $(E,\alpha)$ and $(E',\alpha')$ of elliptic curves over~$k$ equipped with level~$H$ structures are isomorphic when there exists an elliptic curve isomorphism $\phi:E\to E'$ and~$h\in H$ satisfying $\alpha = h^t\circ \alpha'\circ \phi$.

\subsection{The Weil pairing}\label{subsec: Weil pairing}
We also make use of the Weil pairing, as constructed in \cite[§III.8]{Silverman1}. It is a nondegenerate alternating pairing 
$$
e_n:E[n](k)\times E[n](k) \to \mu_n(k),
$$
where $\mu_n(k)$ denotes the group of $n$-th roots of unity in~$k$. If $\alpha: E[n](k)\to (\Z/n\Z)^2$ is a full level~$n$ structure on~$E$ over~$k$, then we define
$$
e_n(\alpha) =  e_n(\alpha^{-1}(1,0)^t, \alpha^{-1}(0,1)^t) \in \mu_n^0(k),
$$
where $\mu_n^0(k)$ denotes the set of primitive $n$-th roots of unity in~$k$.

By \cite[III, Prop.~8.1]{Silverman1}, we have for any $h\in \GL_2(\Z/n\Z)$
$$
e_n(h^t\circ \alpha)^{\det h} = e_n(\alpha).
$$
Consequently, we can associate to a level $H$ structure $\alpha$ the class of $e_n(\alpha)$ in $\mu_n^0(\overline k)/\det H$, where the latter action consists of power-raising. In addition, given an isomorphism $(E,\alpha)\to (E',\alpha')$, we have by \cite[III, Prop.~8.2]{Silverman1}:
$$
e_n(\alpha) = e_n(\alpha') \in \mu_n^0(\overline k)/\det H.
$$
Therefore, $e_n(\alpha)$ is well-defined on isomorphism classes of pairs $(E,\alpha)$.

\subsection{Modular curves}
\label{subsec:modular-fields}

Modular curves arise as coarse moduli spaces of elliptic curves equipped with level structures. This is detailed in Appendix~\ref{appendix}; here we simply record consequences of this characterization over fields.

For any field~$k$ containing~$\Q$, define 
\begin{align*}
    \Pa_H(k) = \{&\text{isomorphism classes of pairs $(E,\alpha)$,}\\
&\text{where $E/k$ and $\alpha$ is a level $H$ structure on $E$}\}.
\end{align*}
The \emph{modular curve} $Y_H$ is a smooth algebraic curve over $\Q$ with the following property: for any field $k$ containing~$\Q$, we have a map
\begin{equation}\label{eq: Pa H k}
\Pa_H(k) \to Y_H(k),
\end{equation}
that is natural with respect to morphisms $k\to k'$, and that is a bijection when $k$ is algebraically closed.

The modular curve $Y_H$ is connected as an algebraic curve over~$\Q$, but is not in general geometrically connected. Its geometric connected components are defined over the number field~$K_H$ constructed as follows. View $\det H$ as a subgroup of~$(\Z/n\Z)^\times$, which we identify to the Galois group of the cyclotomic extension $\Q(\zeta_n)/\Q$. Then~$K_H = \Q(\zeta_n)^{\det H}$ is the fixed field under $\det H$; it is a Galois extension of~$\Q$ with Galois group isomorphic to $(\Z/n\Z)^\times/\det H$, and $\Gal(\Q(\zeta_n)/K_H)\simeq \det H$.

Over~$K_H$, the curve $Y_H$ is a disjoint union of geometrically connected components written as~$Y_H^{(\beta)}$, where $\beta\in (\Z/n\Z)^\times /\det H$, which are smooth algebraic curves over $K_H$. The index $\beta$ corresponds to the different values of the Weil pairing. More precisely, whenever $K_H \subseteq k$, the above map \cref{eq: Pa H k} restricts to a map
$$
\Pa_H^{(\beta)}(k) \to Y_H^{(\beta)}(k),
$$
where $\Pa_H^{(\beta)}(k) \subseteq \Pa_H(k)$ is the set of classes $(E,\alpha)$ satisfying $e_n(\alpha) = \zeta_n^\beta$ in $\mu_n^0(\overline k)/\det H$. Again, this specializes to a bijection when $k$ is algebraically closed. 

The curves $Y_H/\Q$ and $Y_H^{(\beta)}/K_H$ can be compactified into smooth projective curves $X_H/\Q$ and $X_H^{(\beta)}/K_H$ respectively. The complement $X_H^\infty = X_H\backslash Y_H$ is finite, and it is called the set of \emph{cusps} of~$X_H$. The algebraic invariants of (isomorphism classes of) elliptic curves with level structures
that we consider in this paper will be algebraic morphisms of the form $f:X_H^{(-1)}\to \Proj^1$.

\subsection{Modular curves over~\texorpdfstring{$\C$}{C}.}
\label{subsec:complex}

Let~$\UH^* = \UH \cup \mathbb{P}^1(\Q)$ denote the extended upper half plane, on which $\SL_2(\Z)$ acts by Möbius transformations. For a subgroup $H\leq \GL_2(\Z/n\Z)$, define~$\Gamma_H$ as in \eqref{eq:GammaH}. The quotient $\Gamma_H\backslash \UH^*$ is a compact Riemann surface \cite[§2]{Diamond.Shurman2016}, hence an algebraic curve over~$\C$, in which $\Gamma_H\backslash \UH$ is an open dense algebraic subset. We will see that these algebraic curves are exactly the curves $X_H^{(-1)}$ and $Y_H^{(-1)}$, base-changed to~$\C$. By~§\ref{subsec:modular-fields}, $Y_H^{(-1)}(\C)$ is identified with the set of isomorphism classes of elliptic curves over~$\C$ endowed with a level~$H$ structure with associated Weil pairing $\zeta_n^{-1}$; we will give an explicit bijection between $\Gamma_H\backslash \UH$ and this set.

In accordance with \cite[Chapter I]{Silverman2}, for~$\tau\in \UH$, we equip the elliptic curve
$$
E_\tau = \C/(\tau \Z + \Z)
$$
with the ``canonical'' full level~$n$ structure
$$
\begin{matrix}
    \alpha_\tau: &E[n] &\to &(\Z/n\Z)^2,\\[2pt]
    & \frac\tau n &\mapsto &(1,0)^t\\[2pt]
    & \frac 1 n &\mapsto &(0,1)^t.
\end{matrix}
$$
In particular, $e_n(\alpha_\tau)=\zeta_n^{-1}$, by \cite[(2.8.5.3)]{KatzMazur}.

To verify that $\tau\mapsto (E_\tau,\alpha_\tau)$ is the bijection we seek, we can analyze how $\alpha_\tau$ is transformed under  the action of $\SL_2(\Z)$. This will be a particular case of the next two lemmas, which will also be put to use in the next section. 

\begin{lemma}\label{lem: canonical isogeny}
    Let $M\in \GL_2^+(\Q)\cap M_2(\Z)$ be a matrix of determinant $\Delta$ with lower row $(c,d)$, and let $\tau\in \UH$. Then the map
    $$
    \begin{matrix}
        \phi_{M,\tau}: & E_{\tau} &\to &E_{M\tau}\\
        & z &\mapsto &\frac{\Delta}{c\tau + d}z
    \end{matrix}
    $$
    is a well-defined isogeny. Moreover, if $d$ is invertible modulo $\Delta$, then its kernel is given by
    $$
    \ker \phi_{M,\tau} = \left\langle  \frac{c\tau+d}{\Delta}\right\rangle.
    $$
    In particular, $\phi_{M,\tau}$ has degree $\Delta$.
\end{lemma}
\begin{proof}
    Write $M=\begin{pmatrix}a&b\\c&d\end{pmatrix}$, and let
    $\Lambda_\tau=\Z\tau+\Z$ and
    $\Lambda_{M\tau}=\Z(M\tau)+\Z$.
    Since $\Delta=ad-bc$, we have
    $$
    \frac{\Delta}{c\tau+d}\Lambda_\tau = (d(M\tau)-b)\Z + (-c(M\tau)+a)\Z \subseteq \Lambda_{M\tau},
    $$
    proving that $\phi_{M,\tau}$ is well-defined.
    Since multiplication by $\Delta/(c\tau+d)$ is a
    nonzero complex-linear map, it induces a surjective
    holomorphic group homomorphism with finite kernel,
    hence an isogeny.
By definition,
\begin{equation}
[z]\in\ker\phi_M
\iff \frac{\Delta z}{c\tau+d}\in\Lambda_{M\tau}\iff \Delta z\in
\mathbb Z(a\tau+b)+\mathbb Z(c\tau+d).
\label{eq: computation kernel}
\end{equation}
Moreover, as $d$ is invertible modulo $\Delta$, there exists $k\in\Z$
such that
$$
a\equiv kc\mod\Delta,
\qquad
b\equiv kd\mod\Delta,
$$
where the first equivalence comes from $ad \equiv bc \mod \Delta$.
Consequently,
$$
\frac{a\tau+b}{\Delta}
-k\frac{c\tau+d}{\Delta}
\in\Lambda_\tau.
$$
From \cref{eq: computation kernel}, we deduce the result.
\end{proof}
Recall the notation $\pi_n:M_2(\Z)\to M_2(\Z/n\Z)$ for the reduction map modulo~$n$.
\begin{lemma}
\label{lem:transpose-isom}
    Continuing with the notation from \cref{lem: canonical isogeny}, if $\Delta$ is coprime to $n$, then we have
    $$\alpha_{M\tau} \circ \phi_{M,\tau} = (\Delta\pi_n(M)^{-1})^t \circ \alpha_{\tau}.$$
\end{lemma}
\begin{proof}
    This is a direct calculation: writing $M = \begin{pmatrix}
        a & b \\ c & d
    \end{pmatrix}$, we have
    \begin{align*}
        \alpha_{M\tau}  \circ \phi_{M,\tau}(\tau/n) &= \alpha_{M\tau}\left(\frac{\Delta}{c\tau +d} \frac{\tau}{n}\right) \\ &= \alpha_{M\tau}\left(\frac{1}{n} \frac{d(a\tau+b)-b(c\tau+d)}{c\tau +d}\right)\\
        &=\alpha_{M\tau}\left(d\frac{M\tau}{n}-b\frac{1}{n}\right) \\ &=(d,-b)^t 
        \end{align*}
        and
        \begin{align*}
        \alpha_{M\tau}\circ \phi_{M,\tau}(1/n)&=\alpha_{M\tau}\left(\frac{\Delta}{c\tau+d} \frac{1}{n}\right) \\ &= \alpha_{M\tau}\left(\frac{1}{n} \frac{-c(a\tau+b)+a(c\tau+d)}{c\tau+d}\right)\\
        &=\alpha_{M\tau}\left(-c\frac{M\tau}{n}+a\frac{1}{n} \right) \\ &= (-c,a)^t. 
    \end{align*}
    These are the respective images of $\alpha_\tau(\tau/n)$ and $\alpha_\tau(1/n)$ under $(\Delta\pi_n(M)^{-1})^t$.
\end{proof}

The transpose in \cref{lem:transpose-isom} is the reason why we introduced a transposition in the definition of level structures in \cref{def:level-structure}.

Putting~$\Delta=1$ (so that $M\in \SL_2(\Z)$), we obtain as a consequence of \cref{lem:transpose-isom} that $(E_\tau, \alpha_\tau)$ and $(E_{M\tau},\alpha_{M\tau})$ are isomorphic as elliptic curves with a level~$H$ structure if and only if~$\pi_n(M)\in H$, in other words~$M\in \Gamma_H$, as desired. 

We conclude this section by stating that the association $\tau\mapsto(E_\tau,\alpha_\tau)$ indeed stems from an isomorphism of algebraic curves over~$\C$.
\begin{theorem}
\label{thm:complex-modular-curve}
    There exists an isomorphism of complex algebraic curves
    $$
    \Gamma_H\backslash \UH^* \to {(X_H^{(-1)})}_\C
    $$
    under which $\Gamma_H\backslash \mathbb{P}^1(\Q)$ corresponds to the set of cusps of~$X_H^{(-1)}(\C)$, and $\Gamma_H\backslash \UH$ corresponds to ${(Y_H^{(-1)})}_\C$. Moreover, identifying $Y_H^{(-1)}(\C)$ with $\Pa_H^{(-1)}(\C)$ as in \cref{eq: Pa H k}, the above map is given on $\Gamma_H \backslash \UH$ by 
    $$\tau \mapsto (E_\tau, \alpha_\tau).$$
\end{theorem}

As it is not immediate to upgrade the map $\tau\mapsto (E_\tau,\alpha_\tau)$ into an isomorphism of curves, we provide a direct proof building on existing literature, keeping \cref{lem: canonical isogeny,lem:transpose-isom} as motivation and consistency checks only.

\begin{proof}
    Let $G := H\cap \SL_2(\Z/n\Z)$ and $\UH^\pm:=\C \setminus \mathbb{R}$. We have morphisms of complex Riemann surfaces:
    \begin{align*}
    \Gamma_H \backslash \UH &\xrightarrow{\sim} (\UH \times (G \backslash \SL_2(\Z/n\Z)))/\SL_2(\Z) \hookrightarrow (\UH  \times (H \backslash \GL_2(\Z/n\Z))) /  \SL_2(\Z) \\
    & \xrightarrow{\sim} (\UH^\pm \times (H \backslash \GL_2(\Z/n\Z))) / \GL_2(\Z) \xrightarrow{\sim} Y_H(\C).
   \end{align*}
   The first map is defined as in the proof of \cite[Proposition 3.7]{zywina2022possibleindicesgaloisimage} by $\tau \mapsto (\tau, G\cdot I)$; the second is the inclusion as a connected component; the last two are given by \cite[IV. 5.3]{Deligne.Rapoport1973} (after adapting to our conventions). One checks that this composition is indeed given by $\tau \mapsto (E_\tau, \alpha_\tau)$. Since $e_n(\alpha_\tau) = \zeta_n^{-1}$, the image of~$\Gamma_H\backslash \UH$ in $Y_H(\C)$ must be $Y_H^{(-1)}(\C)$. Passing to the compactifications, we see that these morphisms are algebraic by GAGA \cite[§19, Proposition 15]{GAGA}.
\end{proof}

\begin{remark}\label{Gamma and H}
    The information contained in $\Gamma_H$ is coarser than the information contained in $H$; for instance, there may exist two subgroups $H_1 \subsetneq H_2$ of~$\GL_2(\Z/n\Z)$ satisfying $\Gamma_{H_1}=\Gamma_{H_2}$. Then two nonisomorphic $H_1$-structures may become isomorphic as $H_2$-structures, so the natural map $X_{H_1}(\C)\to X_{H_2}(\C)$ is not bijective. It is only thanks to the pairing condition that we have a bijection between $X_{H_1}^{(-1)}(\C)$ and $X_{H_2}^{(-1)}(\C)$ as in \cref{thm:complex-modular-curve}.  
\end{remark}

\section{Modular curves and isogenies}
\label{sec:modular curves and isogenies}

In the absence of level structure, it is well known that $\ell$-isogenies between elliptic curves are parametrized by the modular curve $Y_0(\ell)$, and we have two morphisms $Y_0(\ell) \to Y(1)$ corresponding to the domain and codomain maps in the moduli interpretation. Over~$\C$, this yields a diagram
\begin{equation}
\begin{tikzcd}
    & \Gamma_0(\ell)\backslash \UH \ar[rd, "\tau\mapsto \ell\tau"] \ar[ld, swap, "\tau\mapsto\tau"] & \\
    \SL_2(\Z)\backslash \UH & & \SL_2(\Z)\backslash \UH
\end{tikzcd}
\label{eq:classical-diagram}
\end{equation}
that is the basis for the construction of the classical modular polynomial~$\Phi_\ell$, as in \cite[§11]{Cox1997}. Our aim in this section is to produce a generalization of this diagram when we consider $\ell$-isogenies between elliptic curves endowed with a level~$H$ structure, where~$H\leq \GL_2(\Z/n\Z)$ and~$\ell\nmid n$. This will lay the foundations for the construction of the modular polynomials in §\ref{section: construction}.

In~§\ref{subsection: class isog}, we state that $\ell$-isogenies between elliptic curves with a level~$H$ structure are parametrized by another modular curve $Y_{H'}$ of level dividing $n\ell$. For us, it is crucial to know that the domain and codomain maps $Y_{H'}(\C)\to Y_H(\C)$ indeed arise from algebraic morphisms between modular curves over~$\Q$ (\cref{thm:algebraic morphisms}). We defer the proof to Appendix~\ref{appendix} (\cref{thm:algebraic diag}), as it relies on moduli problems and the associated coarse moduli spaces.

Our second goal, in §\ref{subsec:isog-components}, is to analyze Weil pairings and state that a natural diagram of connected modular curves $Y_{H'}^{(-1)} \rightrightarrows Y_H^{(-1)}$ (as algebraic curves over $K_H$) exists precisely when $(\ell \text{ mod } n)\in \det H$. These maps have degree~$\ell+1$.

After these preparations, we describe the generalization of the diagram~\eqref{eq:classical-diagram} in~§\ref{subsec:complex-isog}, assuming that $(\ell \text{ mod } n)\in \det H$. The domain and codomain morphisms $Y_{H'}^{(-1)}\to Y_H^{(-1)}$, over~$\C$, are also given by $\tau\mapsto\tau$ and $\tau\mapsto M_\ell\tau$ for a suitable matrix $M_\ell\in \GL_2(\Q)^+\cap M_2(\Z)$.

\subsection{Classification of isogenies}\label{subsection: class isog}
Let $E,E'$ be two elliptic curves over a field $k$ of characteristic zero, related by an isogeny $\phi:E\to E'$ of prime degree $\ell\nmid n$. Then~$\phi$ restricts to an isomorphism $E[n]\to E'[n]$. In particular, any level $H$ structure $\alpha':E'[n](\overline k)\to (\Z/n\Z)^2$ on~$E'$ induces unique a level $H$ structure $\phi^*\alpha'=\alpha' \circ \phi$ on~$E$, and conversely. If two pairs $(E,\alpha)$ and $(E',\alpha')$ are already fixed, we say that an isogeny $\phi:E\to E'$ \textit{respects} or \textit{is compatible with} $\alpha$ and $\alpha'$ if $\phi^*\alpha' = \alpha$.

We are thus interested in the moduli problem
\begin{align*}
\Pa_H^{\ell\text{-isog}}(k) = \{&\text{isomorphism classes of triples }(E,\alpha,\phi): \\
&\text{$E$ is an elliptic curve over $k$, $\alpha$ is a level $H$ structure on $E$}, \\
&\text{and $\phi:E\to E'$ is an $\ell$-isogeny}\}, 
\end{align*}
where two triples $(E_1,\alpha_1,\phi_1),(E_2,\alpha_2,\phi_2)$ are declared isomorphic if there are isomorphisms $\psi:(E_1,\alpha_1)\to (E_2,\alpha_2)$ and $\psi':E\to E'$ such that $\psi'\circ \phi_1 = \phi_2 \circ \psi$.

In order to classify solutions to $\Pa_H^{\ell\text{-isog}}$, we begin by identifying $\ell$-isogenies with their kernels. In turn, such subgroups of the $\ell$-torsion are equivalent to  level $B_0(\ell)$ structures, where
\begin{equation}
B_0(\ell):= \left\{\begin{pmatrix}
    * & * \\ 0 & *
\end{pmatrix}\in \GL_2(\Z/\ell\Z) \right\},
\label{eq:B0}
\end{equation}
via
$$
\alpha \mapsto \langle\alpha^{-1}(0,1)^t\rangle.
$$
We then observe that, as $\ell$-isogenies induce level~$H$ structures uniquely, the data of an $\ell$-isogeny respecting level~$H$ structures is equivalent to the data of both a level~$H$ structure and a level~$B_0(\ell)$ structure on the domain. We show in \cref{prop:field isom moduli pb} below that this is equivalent to a level $H'$ structure, where
\begin{equation}\label{eq: def H'}
H' := \crt_{n,\ell}^{-1}(H\times B_0(\ell))
\end{equation}
and
$$\qquad \crt_{n,\ell}:=(\pi_n,\pi_\ell):\Z/n\ell\Z \to \Z/n\Z \times \Z/\ell\Z.
$$
The bijection is explicitly given as follows: to a representative isomorphism
    \[
    \alpha : E[n\ell](\overline k) \to (\Z/n\ell\Z)^2,
    \]
    we associate the triple 
    \begin{equation}\label{rem: moduli equivalence for fields}
    (E, \alpha_n, \phi_\alpha),\qquad\text{where}\qquad \begin{cases}\alpha_n &= (u^{-1}I)^t\circ \alpha|_{E[n](\overline k)}\\
    \phi_\alpha&:E\to E/\langle \alpha_\ell^{-1}(0,1)^t\rangle\\
    \alpha_\ell&= (\pi_\ell(n)^{-1}I)^t\circ \alpha|_{E[\ell](\overline k)}.\end{cases}
    \end{equation}

\begin{proposition}\label{prop:field isom moduli pb}
For any field $k$ with $\Q\subseteq k$, the map\cref{rem: moduli equivalence for fields} gives a bijection
$$
\Pa_{H'}(k)\xrightarrow{\simeq}\Pa_H^{\ell\mathrm{-isog}}(k) 
$$
which is natural in $k$.
In particular, this induces a map
$$
\Pa_H^{\ell\mathrm{-isog}}(k)\to Y_{H'}(k),
$$
also natural in $k$, which is a bijection when $k$ is algebraically closed.
\end{proposition}

\begin{proof}
    See \cref{thm:algebraic diag}.
\end{proof}

We have two natural maps
\begin{equation}\label{eq:diag pb}
\begin{tikzcd}
	& {\Pa_H^{\ell \mathrm{-isog}}(k)} & \\
	{\Pa_H(k)} && {\Pa_H(k)}
	\arrow["{\text{domain}}"{description}, from=1-2, to=2-1]
	\arrow["{\text{codomain}}"{description}, from=1-2, to=2-3]
\end{tikzcd}
\end{equation}
respectively associating to $(E,\alpha,\phi:E\to E')$ the pairs
$$
(E,\alpha) \qquad \text{and}\qquad (E',\alpha \circ \phi^{-1}).
$$
Crucially, these maps arise from morphisms between the corresponding modular curves. The following result is also proven in \cref{thm:algebraic diag}.
\begin{theorem}\label{thm:algebraic morphisms}
    The set-theoretic maps from \cref{eq:diag pb} arise from algebraic morphisms over $\Q$
    \begin{equation}
        \begin{tikzcd}
            & Y_{H'} \ar[ld] \ar[rd] &\\
            Y_H && Y_H.
        \end{tikzcd}
     \label{eq:diagram-YH}
    \end{equation}
    More precisely, for any field $k$ containing $\Q$, the above morphisms induce the following commutative diagram
    \[\begin{tikzcd}
	& {\Pa_H^{\ell\mathrm{-isog}}(k)} & \\
	{\Pa_H(k)} & {Y_{H'}(k)} & {\Pa_H(k)} \\
	{Y_H(k)} && {Y_H(k)}
	\arrow[from=1-2, to=2-1]
	\arrow[from=1-2, to=2-2]
	\arrow[from=1-2, to=2-3]
	\arrow[from=2-1, to=3-1]
	\arrow[from=2-2, to=3-1]
	\arrow[from=2-2, to=3-3]
	\arrow[from=2-3, to=3-3]
\end{tikzcd}\]
    where the map $\Pa_H^{\ell\mathrm{-isog}}(k)\to Y_{H'}(k)$ is from \cref{prop:field isom moduli pb} and $\Pa_H(k)\to Y_H(k)$ is from \cref{eq: Pa H k}.
\end{theorem}
\subsection{Connected components and degrees}
\label{subsec:isog-components}

Since the determinant map $B_0(\ell) \to (\Z/\ell\Z)^\times$ is surjective, the reduction map $\pi_n$ induces an isomorphism
$$
(\Z/n\ell\Z)^\times /\det H' \simeq (\Z/n\Z)^\times/\det H.
$$
Consequently, $K_{H'} = K_H$. Over this field, each of the curves in the diagram~\eqref{eq:diagram-YH} splits as a disjoint union of $\#((\Z/n\Z)^\times/\det H)$ connected components that are also geometrically connected. These components are related as follows.

\begin{proposition}
    For any $\beta\in (\Z/n\Z)^\times/\det H$, the image of~$Y_{H'}^{(\beta)}$ under the maps in~\eqref{eq:diagram-YH} is contained in $Y_H^{(\beta)}\times Y_H^{(\ell \beta)}$.
\end{proposition}

\begin{proof}
    It is sufficient to test this property on $k$-points, where~$k$ is an algebraically closed field containing $K_H$. Let~$\alpha$ be a level $H'$ structure on an elliptic curve~$E/k$; as above, we consider~$\alpha$ as an isomorphism $E[n\ell](\overline k) \to (\Z/n\ell\Z)^2$, well-defined up to the action of $H'$. Further assume that $e_n(\alpha) = \zeta_n^\beta$. Let~$(E',\alpha')$ be the codomain constructed using the diagram~\eqref{eq:diag pb}. Then, by the compatibility between the Weil pairing and isogenies \cite[Chapter III, proposition 8.2.]{Silverman1}, we have
    \[
    e_n(\alpha') = e_n(\alpha)^\ell,
    \]
    so $(E',\alpha')$ corresponds to a point on $Y_H^{(\ell \beta)}(k)$.
\end{proof}

Because our construction of modular polynomials in~§\ref{section: construction} relies on a single invariant on $Y_H^{(-1)}$, we need to assume that $Y_H^{(-1)}=Y_H^{(-\ell)}$. Thus, we only consider subgroups $H\leq \GL_2(\Z/n\Z)$ such that $(\ell \mod n) \in \det H$.
Diagram \cref{eq:diagram-YH} then restricts to a diagram of domain and codomain maps between geometrically connected components of modular curves
\begin{equation}
    \label{eq: diagram connected components}
\begin{tikzcd}
            & Y_{H'}^{(-1)} \ar[ld] \ar[rd] &\\
            Y_H^{(-1)} && Y_H^{(-1)}.
        \end{tikzcd}
\end{equation}
We next analyze the degrees of these maps.
\begin{proposition}\label{prop: degree dom cod}
  Assume $(\ell\mod n)\in H$. Then the domain and codomain morphisms $Y_{H'}^{(-1)} \rightrightarrows Y_H^{(-1)}$ from~\eqref{eq: diagram connected components} have degree $\ell+1$. In particular, the induced maps $Y_{H'}^{(-1)}(\C) \rightrightarrows Y_H^{(-1)}(\C)$ are such that a generic fiber contains $\ell+1$ points.
\end{proposition}
\begin{proof}
    See \cite[Lemma 5.1]{Bröker_2011}. This comes from the fact that level structures are uniquely transported along isogenies. Consequently, fixing either the domain or the codomain, the number of isogenies compatible with a given level structure is the same as the number of isogenies without any level-structure constraints. For the implication on $\C$-points, we apply \cite[Proposition 7.16]{Harris1992}.
\end{proof}
\subsection{Classification of complex isogenies}
\label{subsec:complex-isog}

We now specialize to $k = \C$. Our aim is to construct an analytic uniformization of the diagram~\eqref{eq: diagram connected components}. In other words, we wish to construct the dashed arrows so that the diagram 
\begin{equation}\label{diagram over C}
\begin{tikzcd}
	& {\Gamma_{H'}\backslash\UH} & \\
	{\Gamma_H\backslash\UH} & {Y_{H'}^{(-1)}(\C)} & {\Gamma_H\backslash\UH} \\
	{Y_H^{(-1)}(\C)} && {Y_H^{(-1)}(\C)}
	\arrow[dashed, from=1-2, to=2-1]
	\arrow["\simeq"', from=1-2, to=2-2]
	\arrow[dashed, from=1-2, to=2-3]
	\arrow["\simeq"', from=2-1, to=3-1]
	\arrow["\domain"{description}, from=2-2, to=3-1]
	\arrow["\codomain"{description}, from=2-2, to=3-3]
	\arrow["\simeq"', from=2-3, to=3-3]
\end{tikzcd}
\end{equation}
commutes, where the vertical maps are those from \cref{thm:complex-modular-curve} and the lower ones are those from \cref{eq: diagram connected components} (base-changed to~$\C$).

We will further impose that the dashed arrows are of the form $\tau\mapsto \tau$ and $\tau\mapsto M\tau$ for some matrix~$M\in \GL_2^+(\Q)\cap M_2(\Z)$, and that the $\ell$-isogeny specified by the level~$H'$ structure is precisely the isogeny $\phi_{M,\tau}$ from \cref{lem: canonical isogeny}. The next lemma gathers conditions on~$M$ for this to work.

\begin{lemma}\label{lem:conditions diag}
     Let $\gamma \in \Gamma_0(\ell)$,
    and $M := D(\ell,1)\gamma$.
    For any $\tau \in \UH$, the isogeny $\phi_{M,\tau}:E_\tau \to E_{M\tau}$ from \cref{lem: canonical isogeny} satisfies 
    $$\ker \phi_{M,\tau}=\left\langle \frac 1 \ell \right\rangle.$$ 
   Moreover, $\phi_{M,\tau}$ respects the level~$H$ structures $\alpha_\tau$ and $\alpha_{M\tau}$ if and only if
   \begin{equation}
    \pi_n(\gamma) = uD(u,1)^{-1}M_u^{-1}
    \label{eq: gamma mod n}
    \end{equation}
    for some $M_u\in H$ satisfying $\det M_u = u$.
\end{lemma}
\noindent Note that since $n$ and $\ell$ are coprime and
    $$
    \pi_{n\ell}:\SL_2(\Z) \to \SL_2(\Z/n\ell) \simeq \SL_2(\Z/n\Z) \times \SL_2(\Z/\ell\Z)
    $$
    is surjective, the matrix $\gamma$ can satisfy both $\gamma\in \Gamma_0(\ell)$ and \cref{eq: gamma mod n} simultaneously.
\begin{proof}
   We observe that $\det M=\ell$. In addition, the lower row $(c,d)$ of $M$ is such that $c\equiv 0 \mod \ell$ and $d$ is invertible modulo $\ell$; by \cref{lem: canonical isogeny}, this proves the assertion on $\ker \phi_{M,\tau}$. By \cref{lem:transpose-isom}, $\phi_{M,\tau}$ respects the level structures if and only if $u\pi_n(M)^{-1}\in H$, which is equivalent to \cref{eq: gamma mod n}.
\end{proof}

From now on, we fix a matrix 
$M_u\in H$ such that $\det M_u=u$, and a matrix 
$\gamma_u \in \Gamma_0(\ell)$ such that $\pi_n(\gamma_u) = u D(u,1)^{-1}M_u^{-1}.$
We further define 
\begin{equation}
    \label{def: Ml}
M_\ell:=D(\ell,1)\gamma_u.
\end{equation}
\begin{proposition}\label{prop: maps over C}
    Define $M_\ell$ as in \cref{def: Ml}. Then the map $\Gamma_{H'}\backslash \UH \to Y_{H'}^{(-1)}(\C) \simeq \Pa_H^{\ell\text{-isog}}(\C)$ from \cref{diagram over C} is given by
    $$
    \tau \mapsto \left(\phi_{M_\ell, \tau}:(E_\tau,\alpha_\tau)\to (E_{M_\ell\tau},\alpha_{M_\ell \tau})\right).
    $$
    In particular, any isogeny in $Y_{H'}^{(-1)}(\C)$ is isomorphic to $\phi_{M_\ell,\tau}$ for some $\tau \in \Gamma_{H'}\backslash \UH$.
    Moreover, the domain and codomain maps $Y_{H'}^{(-1)}(\C) \rightrightarrows Y_H^{(-1)}(\C)$ then correspond to the maps $\Gamma_{H'}\backslash \UH \rightrightarrows \Gamma_{H}\backslash \UH$ given by
        $$
        \tau \mapsto \tau \qquad \text{and} \qquad \tau \mapsto M_\ell \tau,
        $$
        in the sense that, defining the upper maps as these, diagram \cref{diagram over C} commutes.
    \end{proposition}
\begin{proof}
By \cref{rem: moduli equivalence for fields}, the image of $\tau$ under $\Gamma_{H'}\backslash \UH \to Y_{H'}^{(-1)}(\C)$ is equal to the quotient isogeny 
    $$
    (E_\tau,\alpha_\tau) \xrightarrow{q} (E_\tau/ C, \alpha_\tau\circ q^{-1}), \qquad
    \text{ where } \qquad C = \langle \alpha_{\tau,\ell}^{-1}(0,1)^t\rangle = \left\langle \frac{1}{\ell}\right\rangle.$$
    Note that the latter equality is the reason why we chose this normalization in \cref{rem: moduli equivalence for fields}.
    By \cref{lem:conditions diag}, the isogeny $\phi_{M_\ell,\tau}:(E_\tau,\alpha_\tau)\to(E_{M_\ell\tau},\alpha_{M_\ell\tau})$ is compatible with the level~$H$ structures. In addition, this lemma also states that $\phi_{M,\tau}$ and $q$ have the same kernel, hence are isomorphic. Since they both transport level structures to $\alpha_\tau$, they are also isomorphic as structure-preserving isogenies.
\end{proof}
Paving the way for §\ref{subsec: construction}, we can now use the above result to encode all possible codomains of isogenies that have a fixed domain in $X_H^{(-1)}(\C)$.
Since $X_{H'}^{(-1)}$ and $X_H^{(-1)}$ are smooth projective curves, the degree-$(\ell+1)$ morphism 
$$\domain:X_{H'}^{(-1)}\to X_H^{(-1)}$$ 
induces a morphism
$$
\domain^{-1}:X_{H}^{(-1)} \to \Sym^{\ell+1}\left(X_{H'}^{(-1)}\right),
$$
where $\Sym^{\ell+1}$ denotes the $(\ell+1)$-th symmetric power, assigning a point to its fiber. In the moduli interpretation, this amounts to considering all $\ell$-isogenies with a given domain. Further, composing with the codomain morphism $X_{H'}^{(-1)}\to X_H^{(-1)}$, we get a morphism
\begin{equation}
\codomain\circ \domain^{-1}:X_H^{(-1)} \to \Sym^{\ell+1}\left(X_H^{(-1)}\right).
\label{eq: cod dom -1}
\end{equation}

\begin{corollary}\label{cor: domains and codomains}
    The morphism \cref{eq: cod dom -1} is given on $\C$-points by the map
    $$
    \Gamma_H \backslash \UH^* \to \Sym^{\ell+1}(\Gamma_H \backslash \UH^*),\qquad\tau \mapsto (\gamma\tau)_{\gamma \in \Gamma_H \backslash (\Gamma_H M_\ell \Gamma_H)}.
    $$
    In particular, the above map depends neither on the choice of $M_u$ or $\gamma_u$. Moreover, for generic $\tau$, the $\ell+1$ elements of the above family are distinct in $\Gamma_H \backslash \UH^*$.
\end{corollary}
\begin{proof}

For a fixed $\tau\in \Gamma_H \backslash \UH^* $, by \cref{prop: maps over C}, the isogenies with domain $\tau$ are parametrized by the classes in $\Gamma_{H'}\backslash \UH^*$ given by
$$
(\gamma\tau)_{ \gamma \in  \Gamma_{H'}\backslash\Gamma_H}.
$$
It then suffices to use the bijection
$$
\Gamma_{H'}\backslash\Gamma_H
\xrightarrow{\sim}
\Gamma_H\backslash(\Gamma_HM_\ell\Gamma_H),
\qquad
\Gamma_{H'}\gamma
\longmapsto
\Gamma_HM_\ell\gamma.
$$
Concerning multiplicities, since, by construction, $-I\in \Gamma_H$ if and only $-I\in \Gamma_{H'}$, there can be no distinct representatives $\gamma,\gamma'\in \Gamma_H\backslash\Gamma_H M_\ell \Gamma_H$ satisfying $\gamma\tau = \gamma'\tau$ for every $\tau\in \UH$. Hence the codomains are pairwise non-isomorphic for generic~$\tau$.
\end{proof}

\section{Modular polynomials for level structures}
\label{section: construction}

In this section, we carry out the construction and study of modular polynomials, building on the modular curves framework of §\ref{section: modular curves}--\ref{sec:modular curves and isogenies}. Throughout, we retain a fixed subgroup $H\leq \GL_2(\Z/n\Z)$ satisfying $(\ell\mod n)\in \det H$ for a fixed prime $\ell\nmid n$, and add the condition that $X_H^{(-1)}$ has genus 0. We then fix an isomorphism of curves $f:X_H^{(-1)}\to \Proj_{K_H}^1$ such that $f(\infty)=\infty$ over $\C$. We also denote by $f$ the induced modular function 
$$
\UH \to \Gamma_H \backslash \UH \to Y_H^{(-1)}(\C) \to \Proj^1(\C).
$$
In §\ref{subsec: construction}, we construct the modular polynomial $\Phi_\ell^f$, proving \cref{thm:main-existence-1,thm:main-existence-2}. In §\ref{subsec: action cusps}, we then determine a simple and equivalent condition on $H$ for $\Phi_\ell^f$ to be denominator-free, proving  \cref{thm: summary conditions} \ref{item: intro denominator}. Assuming that this condition holds, we study the symmetry of $\Phi_\ell^f$ in~§\ref{subsec: symmetry}; then we find explicit upper-triangular representatives of the cosets $\Gamma_H\backslash \Gamma_H M_\ell \Gamma_H$ appearing in the expression of $\Phi_\ell^f$ in §\ref{subsec: explicit representatives}, and deduce sufficient conditions to ensure that $\Phi_\ell^f$ has integer coefficients.

\subsection{Construction of the modular polynomial}\label{subsec: construction}

For $\tau \in \UH$, consider the intermediary polynomial
\begin{equation}\label{def Phi^tau}
\Phi_{\ell}^{f,\tau} (y)=\prod_{\gamma \in \Gamma_H \backslash (\Gamma_H M_\ell \Gamma_H)} (y-f(\gamma \tau)),
\end{equation}
which is well defined since $f^{-1}(\infty)=\infty$ is a cusp.
By \cref{cor: domains and codomains}, the univariate polynomial $\Phi_\ell^{f,\tau}$ does not depend on the choices of $M_u$ or $\gamma_u$.
Write 
$$\Phi_\ell^{f,\tau}(y)=\sum_{r=0}^{\ell+1}C_r(\tau) y^r.$$
\begin{lemma}
    The coefficients $C_r(\tau)$ are modular functions for $\Gamma_H$ that are holomorphic on $\UH$. In particular, they are rational fractions in $f$.
\end{lemma}
\begin{proof}
Let $\gamma\in\Gamma_H$. Right multiplication by $\gamma$ gives a bijection
$
\Gamma_H\backslash(\Gamma_HM_\ell\Gamma_H)
\to
\Gamma_H\backslash(\Gamma_HM_\ell\Gamma_H),
$
whose inverse is right multiplication by $\gamma^{-1}$. Hence we have
$$
\Phi_\ell^{f,\gamma\tau}(y)
=
\Phi_\ell^{f,\tau}(y).
$$
Comparing the coefficients of $y^r$ gives
$C_r(\gamma\tau)=C_r(\tau)$, so $C_r$ is $\Gamma_H$-invariant.
Moreover, each $C_r$ is a polynomial in the functions
$\tau\mapsto f(\gamma_*\tau)$ appearing above. Since $f$ is holomorphic on
$\UH$, so is each $f(\gamma_*\tau)$, and therefore so is $C_r$.
Likewise, since $f$ is meromorphic at the cusps, each $f(\gamma_*\tau)$ and hence
each $C_r$ is meromorphic at the cusps. Thus $C_r$ is a modular function for
$\Gamma_H$ that is holomorphic on $\UH$. We conclude since $f$ is a Hauptmodul (see for example the argument in \cite[chapter I, corollary 4.2 (a)]{Silverman2}).
\end{proof}

\begin{definition}
For each $0\leq r \leq \ell+1$, let $R_r(x)\in \C(x)$ be the rational fraction such that 
$$
C_r(\tau) = R_r(f(\tau)).
$$
The rational fraction
\begin{equation}
\Phi_\ell^f(x,y):=\sum_{r=0}^{\ell+1}R_r(x)y^r \in \C(x)[y]
\label{eq: def Phi ell f}
\end{equation}
is called the \textit{modular polynomial of invariant $f$ and degree $\ell$}.
Putting all fractions over a common denominator, there are polynomials $A(x,y)\in \C[x,y]$, $B(x)\in \C[x]$ such that
$$
\Phi_\ell^f(x,y) = \frac{A(x,y)}{B(x)}.
$$
We say that $\Phi_\ell^f$ is \emph{denominator-free} if $\deg B=0$.
\end{definition}

The next theorem summarizes fundamental properties of~$\Phi_\ell^f$, and includes \cref{thm:main-existence-1,thm:main-existence-2}.

\begin{theorem}\label{thm:properties-Phi-ell-f}
    The rational fraction $\Phi_\ell^f$ has the following properties:
    \begin{enumerate}
    \item \label{item: Phi=Phi tau}
    For any $\tau\in \UH$, we have 
    \begin{equation}\label{Phi=Phi tau}
\Phi_\ell^f(f(\tau),y)=\Phi_{\ell}^{f,\tau} (y)=\prod_{\gamma \in \Gamma_H \backslash (\Gamma_H M_\ell \Gamma_H)} (y-f(\gamma \tau)).
\end{equation}
    \item \label{item: classification Phi} For any $\tau,\tau'\in \UH$,  there exists an $\ell$-isogeny $(E_\tau,\alpha_\tau)\to (E_{\tau'},\alpha_{\tau'})$ if and only if 
    \[\Phi_\ell^f(f(\tau),f(\tau'))=0.\]
    In particular, for any two pairs $(E,\alpha),(E',\alpha')$ with pairing $\zeta_n^{-1}$, there exists an $\ell$-isogeny $(E,\alpha)\to(E',\alpha')$ if and only if there exist $\tau,\tau'\in \UH$ such that
        $(E,\alpha)\simeq (E_\tau,\alpha_\tau)$, $(E',\alpha')\simeq (E_{\tau'},\alpha_{\tau'})$, and $\Phi_\ell^f(f(\tau),f(\tau'))=0$.
    \item \label{item: hypothesis 1.1}
    Any two elliptic curves~$E,E'$ over~$\C$ are $\ell$-isogenous if and only if
  there exist $\tau,\tau'\in \UH$ satisfying
  $$E\simeq E_\tau,\quad E'\simeq E_{\tau'}, \quad\text{and}\quad \Phi_\ell^f(f(\tau),f(\tau')) = 0.$$
    
    \item \label{item: degree of Phi_ell^f} We have $\deg_y\Phi_\ell^f=\ell+1$ and $\deg_xA(x,y)=\ell+1.$
    In addition, the coefficient of $y^{\ell+1}$ in $\Phi_\ell^f$ is 1.
    \item \label{item: equivalence property star} The rational fraction $\Phi_\ell^f$ is denominator-free if and only if
    \begin{numberedcenter}\label{property star}
    $M_\ell^{-1} \Gamma_H \infty \subseteq \Gamma_H \infty$.
    \end{numberedcenter}
    In particular, if $X_H^{(-1)}$ has a single cusp, then $\Phi_\ell^f$ is denominator-free.
     \item \label{item: coefs Phi} The rational fraction $\Phi_\ell^f$ has coefficients in $K_H$.
     \item \label{item: B(y)} The polynomial $B(x)$ may only vanish at images of cusps under $f$.
     \item \label{item: Phi irreducible} If $\Phi_\ell^f$ is denominator-free, then it is an irreducible polynomial over $\C$.
    \end{enumerate}
\end{theorem}
\begin{proof}

We prove each item in order.
\begin{enumerate}
\item Equation \cref{Phi=Phi tau} holds by construction.

\item Given $\tau,\tau'\in \UH$, by \ref{item: Phi=Phi tau}, $\Phi_\ell^f(f(\tau),f(\tau'))=0$, if and only if $f(\tau')=f(\gamma\tau)$ for some $\gamma \in \Gamma_H \backslash \Gamma_H M_\ell \Gamma_H$. This holds if and only if $(E_{\tau'},\alpha_{\tau'})\simeq(E_{\gamma \tau},\alpha_{\gamma \tau})$, which is equivalent to the existence of an $\ell$-isogeny $(E_\tau,\alpha_\tau)\to (E_{\tau'},\alpha_{\tau'})$.
In particular, given two elliptic curves $(E,\alpha),(E',\alpha')$ over $\C$ with level $H$ structure of pairing $\zeta^{-1}$, if $\tau,\tau'$ exist with the desired properties, then there exists an isogeny $(E,\alpha)\to (E',\alpha')$. Conversely, by \cref{prop: maps over C}, any isogeny $(E,\alpha)\to (E',\alpha')$ is isomorphic to $\phi_{M_\ell,\tau}:(E_\tau,\alpha_\tau)\to (E_{M_\ell\tau},\alpha_{M_\ell \tau})$ for some $\tau\in \Gamma_{H'}\backslash\UH$, implying that $\tau,M_\ell \tau$ have the desired properties.

\item 
 Let $E,E'$ be two elliptic curves over $\C$. Then the following assertions are equivalent:
    \begin{enumerate}
        \item there exists an $\ell$-isogeny $E\to E'$,
        \item there exist level $H$ structures $\alpha,\alpha'$ and an $\ell$-isogeny $(E,\alpha)\to (E',\alpha')$,
        \item there exist structures $\alpha,\alpha'$ and $\tau,\tau'\in \UH$ such that $(E,\alpha)\simeq (E_\tau,\alpha_\tau),$ $(E',\alpha')\simeq (E_{\tau'},\alpha_{\tau'})$ and $\Phi_\ell^f(f(\tau),f(\tau'))=0$,
        \item there are $\tau,\tau'\in \UH$ so that $E\simeq E_\tau$, $E'\simeq E_{\tau'}$  and $\Phi_\ell^f(f(\tau),f(\tau'))=0$.
    \end{enumerate}
    The first and the third equivalences hold by inducing level structures, and the second by item~\ref{item: classification Phi}.

\item By construction, from \cref{Phi=Phi tau}, we have 
$$\deg_y \Phi_\ell^f = \deg \Phi_\ell^{f,\tau}=\ell+1,$$ 
where $\tau\in \UH$ is arbitrary. Moreover, fixing a generic $\tau'\in \Gamma_H \backslash \UH$, by \ref{item: classification Phi}, the numbers $x\in \C$ satisfying $A(x,f(\tau'))=0$ are given by the elements $f(\tau)$, where $\tau\in \Gamma_H \backslash \UH$ is such that $(E_{\tau},\alpha_{\tau})$ is the domain of an isogeny with codomain $(E_{\tau'},\alpha_{\tau'})$. In other words, the elements $\tau$ correspond to the domains of the isogenies mapping to $\tau'$ under the codomain morphism $$\Gamma_{H'} \backslash \UH\simeq Y_{H'}^{(-1)}(\C) \to \Gamma_H\backslash \UH \simeq Y_H^{(-1)}(\C).$$
By \cref{prop: degree dom cod}, there are $\ell+1$ such preimages, proving that 
$$
\deg_x A(x,y)=\ell+1.$$
In addition, by construction, from \cref{eq: def Phi ell f}, the coefficient of $y^{\ell+1}$ in $\Phi_\ell^f$ is 1.
\item We have $\deg B(x)=0$, if and only if $R_r(x)$ is denominator-free for every $1\leq r \leq \ell+1$. Recall that $R_r$ satisfies $C_r(\tau)=R_r(f(\tau))$. 
For any $\gamma \in \Gamma_H \backslash (\Gamma_H M_\ell \Gamma_H)$, let $f_\gamma(\tau):=f(\gamma \tau)$ and consider $C_r$ and every $f_\gamma$ as analytic functions $\UH^* \to \Proj^1(\C)$. In particular, as $f$ is holomorphic on $\UH$, $f_\gamma$ is also holomorphic on $\UH$.
By the same argument as in \cite[chapter I, corollary 4.2 (b)]{Silverman2}, $R_r$ is denominator-free if and only if $C_r(\tau)$ is holomorphic on $\UH^* \setminus (\Gamma_H \infty)$. In turn, this is equivalent to each $f_\gamma$ having no pole outside of $\Gamma_H \infty$. In addition, as $f_\gamma$ is holomorphic on $\UH$, this is equivalent to each $f_\gamma$ having a pole at some $\kappa \in \Proj^1(\Q)$ only if $\kappa\in \Gamma_H \infty$. Moreover $f_\gamma$ has a pole at $\kappa$ if and only if $f$ has a pole at $\gamma \kappa$; equivalently, $\gamma\kappa \in \Gamma_H\infty$. Therefore, $R_r$ is denominator-free if and only if:
\begin{numberedcenter}\label{ppte intermediaire}
for all $\gamma\in \Gamma_H \backslash (\Gamma_H M_\ell \Gamma_H)$, any $\kappa\in \Proj^1(\Q)$\\ such that $\gamma\kappa \in \Gamma_H \infty$ satisfies $\kappa \in \Gamma_H\infty$. 
\end{numberedcenter}
\par\noindent In particular, if $X_H^{(-1)}$ admits a single cusp, then $\Proj^1(\Q) = \Gamma_H \infty$, so the above property holds.
Writing $\gamma=M_\ell \gamma'$ for $\gamma'\in \Gamma_H$ and $\kappa=\gamma'^{-1}\kappa'$, \cref{ppte intermediaire} is equivalent to having:
\begin{center}
for all $\gamma'\in \Gamma_H$, any $\kappa'\in \Proj^1(\Q)$ such that $M_\ell \kappa' \in \Gamma_H\infty$ satisfies $\gamma'^{-1}\kappa' \in \Gamma_H\infty$ (or equivalently, $\kappa'\in \Gamma_H\infty$).
\end{center}
Ultimately, writing $M_\ell\kappa' = \gamma''\infty$ for $\gamma''\in \Gamma_H$, this equivalent to having that for any $\gamma''\in \Gamma_H$, there exists $\gamma_*\in \Gamma_H$ such that $M_\ell^{-1}\gamma''\infty = \gamma_* \infty$.

\item Recall that $f$ and the domain and codomain morphisms $X_{H'}^{(-1)}\rightrightarrows X_H^{(-1)}$ are defined over $K_H$. Consequently, the image scheme of the morphism 
\begin{equation}
X_{H'}^{(-1)} \to X_H^{(-1)} \times X_H^{(-1)}\xrightarrow{f\times f}\Proj^1 \times \Proj^1
\label{eq: image variety}
\end{equation}
is also defined over $K_H$.
In addition by \ref{item: classification Phi}, the $\C$-points of this image scheme form the locus of $\Phi_\ell^f$. Moreover, $\Phi_\ell^f$ has leading coefficient 1 in $y$. Therefore, it has coefficients in $K_H$.

\item If $B$ vanishes at $x$, then for some $1\leq r\leq \ell+1$, $R_r$ has a pole at $x$. Therefore, $C_r(\tau)$ has a pole at $f^{-1}(x)$, implying that $f^{-1}(x)$ is a cusp.
\item Since the locus of $\Phi_\ell^f$ is equal to the image of the morphism \cref{eq: image variety} and $X_{H'}^{(-1)}$ is irreducible, $\Phi_\ell^f$ must be of the form $c\Phi'(x,y)^m$ for some $c\in K_H$, $\Phi'(x,y)\in K_H[x,y]$ irreducible and $m\geq 1$. In addition, by \cref{cor: domains and codomains} and \ref{item: Phi=Phi tau}, for generic $\tau \in \UH$, 
$$
c\Phi'(f(\tau),y)^m=\Phi_\ell^f(f(\tau),y)=\Phi_\ell^{f,\tau}(y)
$$
is squarefree over $\C$, meaning that $m=1$ and $\Phi_\ell^f$ is irreducible over $\C$. \qedhere
\end{enumerate}
\end{proof}

In particular, $\Phi_\ell^f$ satisfies the hypotheses of \cref{thm:main-height} if \cref{property star} holds.

\subsection{Action on cusps}\label{subsec: action cusps}

In this section, we prove \cref{thm: summary conditions} \ref{item: intro denominator}, giving a necessary and sufficient condition for the polynomial $\Phi_\ell^f$ to be denominator-free.
We translate the arithmetic condition  \cref{property star} on the $\Gamma_H$-cusps into a group-theoretic condition on $H$ and $u=(\ell \mod n)$. 

Let $e_1=(1,0)^t$, either in $\Z^2$ or in $(\Z/n\Z)^2$. For a matrix $M$, we sometimes denote by $M_1=Me_1$ its first column.
For any coprime integers $a,c$, let
$$
\delta_\ell(a,c)^t:= \begin{cases}
    (a,\ell c)^t & \text{ if } \ell \nmid a\\
    (a/\ell,c)^t& \text{ if } \ell \mid a.
\end{cases}
$$
Recall that we take $\gamma_u\in \Gamma_0(\ell)$ satisfying $\pi_n(\gamma_u) = uD(u,1)^{-1}M_u^{-1}$, where $M_u\in H$ has determinant $u$, and $M_\ell = D(\ell,1)\gamma_u
$. 

\begin{lemma}\label{lem: 1st reformulation}
The following conditions are equivalent:
\begin{enumerate}
    \item \label{condition cusps} We have $M_\ell^{-1} \Gamma_H \infty \subseteq \Gamma_H \infty$,
    \item \label{upstairs upstairs}
     For any $\gamma \in \Gamma_H$, there exists~$\varepsilon\in \{\pm 1\}$ such that $\varepsilon \gamma_u^{-1}\delta_{\ell}(\gamma_1)\in \Gamma_He_1$.
\end{enumerate}
\end{lemma}
\begin{proof}
 For a vector of coprime integers $(a,c)^t\in \Z^2$, denote by $(a,c)^t \infty:=a/c\in \Proj^1(\Q)$. 
 Therefore, for any $\gamma \in \Gamma_H$, we have 
 $$M_\ell^{-1}\gamma \infty= (M_\ell^{-1}\gamma)_1 \infty=(M_\ell^{-1}\gamma_1) \infty.$$
Write 
$$\gamma_u = \begin{pmatrix}
    a_0 & b_0 \\ c_0 & d_0
\end{pmatrix},\qquad\gamma_u^{-1}=\begin{pmatrix}
    d_0 & -b_0 \\ -c_0 & a_0
\end{pmatrix},\qquad\gamma=\begin{pmatrix}
    a & b \\ c & d
\end{pmatrix}.$$
We calculate 
$$M_\ell^{-1}\gamma_1=\gamma_u^{-1}D(\ell,1)^{-1}\gamma_1=(d_0a/\ell-b_0c,-c_0a/\ell+a_0c)^t.$$
On the other hand, we have
$$
\gamma_u^{-1}\delta_\ell(a,c)^t
    =\begin{cases}
    (d_0a-\ell b_0c,-c_0a+\ell a_0c)^t & \text{ if } \ell \nmid a\\
    (d_0a/\ell-b_0c,-c_0a/\ell+a_0c)^t & \text{ if } \ell \mid a.
\end{cases}
$$
 Therefore, $$(M_\ell^{-1}\gamma_1) \infty =\frac{d_0a/\ell-b_0c}{-c_0a/\ell+a_0c}= \gamma_u^{-1}\delta_{\ell}(a,c)^t\infty.$$
 Consequently, we can reformulate condition \ref{condition cusps} as
 \begin{numberedcenter}\label{goal}
 for any $\gamma\in \Gamma_H$, there exists $\gamma_*\in \Gamma_H$ such that
     $\gamma_u^{-1}\delta_{\ell}(\gamma)_1 \infty = (\gamma_*)_1\infty$.
\end{numberedcenter}     
Since multiplication by elements in $\SL_2(\Z)$ preserves the $\gcd$, we have
 $$\gcd (\gamma_u^{-1}\delta_{\ell}(a,c)^t)=\gcd \delta_\ell(a,c)^t=1.$$
 Thus, given $\gamma\in \Gamma_H$, some $\gamma_*\in \Gamma_H$ satisfies \cref{goal} if and only if $(\gamma_*)_1 = \pm \gamma_u^{-1}\delta_{\ell}(\gamma)_1$. Therefore, \cref{goal} is equivalent to condition \ref{upstairs upstairs}.
\end{proof}

Condition \ref{upstairs upstairs} from \cref{lem: 1st reformulation} means that for any $\gamma\in \Gamma_H$, there exists $\varepsilon \in \{\pm 1\}$ such that the vector $\varepsilon\gamma_u^{-1}\delta_\ell(\gamma_1)$ is the first column of a matrix in $\Gamma_H$. We now show that finding such a matrix in $\Gamma_H$ can be reduced to finding one in $\pi_n(\Gamma_H)=H\cap \SL_2(\Z/n\Z)$.
\begin{lemma}\label{lemma reduction mod n}
 Let $\Gamma\leq \SL_2(\Z)$ be a subgroup with $\Gamma(n) \subseteq \Gamma$ and $\mathrm{v}\in \Z^2$ be a column vector of coprime numbers. Then $\mathrm{v}\in \Gamma e_1$ if and only if $\pi_n(\mathrm{v}) \in \pi_n(\Gamma)e_1$.
    \end{lemma}
\begin{proof} Assuming that $\mathrm{v}\in \Gamma e_1$, we immediately have $\pi_n(\mathrm{v}) \in \pi_n(\Gamma)e_1$.
    Suppose that $\pi_n(\mathrm{v}) \in \pi_n(\Gamma)e_1$, write 
    $\mathrm{v}=(a,c)^t$ and let $b_0,d_0\in \Z$ be such that 
    $$g=\pi_n \begin{pmatrix}
        a & b_0 \\ c & d_0
    \end{pmatrix}\in \pi_n(\Gamma).$$
    As $a d_0-b_0c\equiv 1\mod n$, there exists $m \in \Z$ such that
    $$
    ad_0-b_0c+nm = 1.
    $$
    As $\gcd(a,c)=1$, let $k_1,k_2\in \Z$ be such that
    $$
    k_2a+k_1c=m.
    $$
    Then, letting 
    $$b=b_0-nk_1, \qquad d = d_0+nk_2, \qquad \gamma_*:=\begin{pmatrix}
        a & b \\ c & d
    \end{pmatrix},$$ 
    the matrix $\gamma_*$ satisfies $\gamma_*e_1=\mathrm{v}$ and
    $$
    \det \gamma_* = ad_0-b_0c+n(k_2a+k_1c)=ad_0-b_0c+nm =1.
    $$
    In addition, it also satisfies $\pi_n(\gamma_*)=g \in \pi_n(\Gamma)$, implying that 
    \begin{equation} \gamma_*\in \pi_n^{-1}\pi_n(\Gamma) \cap \SL_2(\Z) = \Gamma. \tag*{\qedhere}\end{equation}
\end{proof}

For the next proposition, we define
$$U(u,v)=\begin{pmatrix}
    u^{-1} & v \\ 0 & u
\end{pmatrix}\quad\text{and}\quad M(u,w)=\begin{pmatrix}
    1 & w \\ 0 & u
\end{pmatrix}.$$

\begin{proposition}[= \cref{thm: summary conditions} \ref{item: intro denominator}]\label{equivalence}
Keeping notation from \cref{lem: 1st reformulation}, the following conditions are equivalent:
\begin{enumerate}
     \item \label{item: cond d-f in thrm} We have $M_\ell^{-1}\Gamma_H \infty \subseteq \Gamma_H \infty$,
    \item There exist $v,w\in \Z/n\Z$ and $\varepsilon_1,\varepsilon_2\in \{\pm 1\}$ such that
    \begin{equation}
    \varepsilon_1 U(u,v), \ \varepsilon_2 M(u,w) \in H.
    \label{eq: matrix condition ppty star}
    \end{equation}
\end{enumerate}
\end{proposition}
\begin{proof}
 Using \cref{lem: 1st reformulation} and \cref{lemma reduction mod n}, we can reduce condition \ref{item: cond d-f in thrm} modulo $n$.
It is equivalent to:
\begin{numberedcenter}\label{upstairs downstairs}
for any $\gamma\in \Gamma_H$, there exists $\varepsilon\in \{\pm 1\}$ such that $\pi_n(\varepsilon \gamma_u^{-1}\delta_{\ell}(\gamma_1))\in G_He_1$,
\end{numberedcenter}
\par\noindent where $G_H := \pi_n(\Gamma_H)=H\cap \SL_2(\Z/n\Z)$.
We then prove that \cref{eq: matrix condition ppty star} and \cref{upstairs downstairs} are equivalent. 
We compute
\begin{equation}
\pi_n(\varepsilon \gamma_u^{-1}\delta_{\ell}(\gamma_1))=\varepsilon \delta M_u\pi_n(\gamma_   1), \qquad  \text{where }\delta = \begin{cases}
    1 & \text{ if } \ell \nmid \gamma_{1,1}\\
    u^{-1}& \text{ if } \ell \mid \gamma_{1,1}.
\end{cases} 
\label{eq: computation pi delta}
\end{equation}
Assume first \cref{upstairs downstairs}.
Let $\gamma,\gamma'\in \Gamma_H$ be such that $\pi_n(\gamma)= \pi_n(\gamma')=I$ and  $\ell \nmid \gamma_{1,1}$, $\ell \mid \gamma'_{1,1}$. Applying \cref{upstairs downstairs} to these and using \cref{eq: computation pi delta}, there exist $g,g'\in G_H$ and $\varepsilon,\varepsilon' \in \{\pm 1\}$ such that 
$$ge_1 = \varepsilon M_ue_1\qquad g'e_1=\varepsilon' u^{-1}M_ue_1.$$ 
This implies that $g^{-1}g'e_1 = \varepsilon\varepsilon' u^{-1} e_1$.
Since $g^{-1}g' \in \SL_2(\Z/n\Z)$, there exists $v\in \Z/n\Z$ such that
\begin{equation}
g^{-1}g' = \varepsilon \varepsilon' \begin{pmatrix}
    u^{-1} & v \\ 0 & u
\end{pmatrix}\in H.
\label{eq: matrix 1}
\end{equation}
On the other hand, since $M_u^{-1}g e_1 = \varepsilon e_1$ and $\det (M_u^{-1}g) = u^{-1}$, there exists $w'\in \Z/n\Z$ such that 
$$
M_u^{-1}g = \varepsilon \begin{pmatrix}
    1 & w' \\ 0 & u^{-1}
\end{pmatrix}, \qquad\text{hence}\qquad \varepsilon\begin{pmatrix}
    1 & -w'u \\ 0 & u
\end{pmatrix}\in H.
$$
This concludes the first implication.
Assume now \cref{eq: matrix condition ppty star} and let 
$$U=\varepsilon_1\begin{pmatrix}
    u^{-1} & v \\ 0 & u
\end{pmatrix},\qquad M=\varepsilon_2\begin{pmatrix}
    1 & w \\ 0 & u
\end{pmatrix}\in H,$$
where $\varepsilon_1,\varepsilon_2\in \{\pm 1\}$. Since $\det U=1$, we have $U\in G_H$. Therefore,
\begin{equation}
\varepsilon_1u^{-1}e_1 = Ue_1 \in G_H e_1.
\label{eq: ue1}
\end{equation}
In addition, since $M\in H$ and $\det M = u = \det M_u$, we have $g:=M_uM^{-1} \in G_H$. As $\varepsilon_2Me_1 = e_1$, this gives
\begin{equation}\label{eq: ge1}
    \varepsilon_2M_ue_1 = \varepsilon_2gMe_1=ge_1\in G_He_1.
\end{equation}
Therefore, from \cref{eq: ue1} and \cref{eq: ge1}, we have
$$
\varepsilon_1\varepsilon_2u^{-1}M_ue_1= \varepsilon_1u^{-1}ge_1=g\varepsilon_1u^{-1}e_1= gUe_1\in G_He_1.
$$
Let $\gamma \in \Gamma_H$ and $h=\pi_n(\gamma)\in G_H$. Let $h' = M_uhM_u^{-1}$, which satisfies $h'\in G_H$.
The vectors
$$
\varepsilon_2M_uh_1 = h'(\varepsilon_2 M_ue_1) =h'ge_1,
\qquad
\varepsilon_1\varepsilon_2 u^{-1}M_uh_1 = h'(\varepsilon_1\varepsilon_2 u^{-1}M_ue_1)=h'gUe_1
$$
then belong to $G_He_1$.
Therefore, whichever case in \cref{eq: computation pi delta} occurs, we have $$\pi_n(\varepsilon \gamma_u^{-1}\delta_\ell(\gamma_1))\in G_He_1$$
for some $\varepsilon \in \{\pm 1\}$.
\end{proof}

\begin{example}
    In the classical case, $n=1$, so $$\SL_2(\Z/n\Z)=\GL_2(\Z/n\Z)=H=\{I\}.$$ Also, $(\ell \mod n)=u=1\in \det H$ for any $\ell$, so \cref{eq: matrix condition ppty star} is satisfied. Therefore, the classical modular polynomial $\Phi_\ell=\Phi_\ell^j$ is denominator-free.
    
    As a second example, let $H =B_0(n)$. Then $(\Z/n\Z)^\times \subseteq \det H$. We have $U(u,0)\in H$ and $D(1,u)=M(u,0)\in H$, meaning that any modular polynomial constructed from an invariant on $X_H^{(-1)}$ would be denominator-free.
    
    For a non-example, take  $n=5$, $\ell=2$, and again $M_u=D(1,u)$ and $H=\langle M_u\rangle$. Then $H\cap \SL_2(\Z/n\Z)=\{I\}$ and $\Gamma_H = \Gamma(5)$. However, $I$ is not of the form $\pm U(u,v)$. Therefore, in that case \cref{eq: matrix condition ppty star} is not satisfied and by \cref{thm:properties-Phi-ell-f} \ref{item: equivalence property star}, for no isomorphism $f:X_{H}^{(-1)} \xrightarrow{\sim} \Proj^1_{K_H}$, is $\Phi_\ell^f$ denominator-free.
\end{example}
\begin{remark}\label{cor: setting denominator free}
    If $\Phi_\ell^f$ is denominator-free, then as $\Phi_\ell^f$ does not depend on the choices of $M_u$ and $\gamma_u$, we can make the following advantageous choices, in line with \cref{prop: maps over C}:  $$M_u=\varepsilon_2M(u,w)\quad\text{and}\quad\gamma_u=\varepsilon_2T^{\chi},$$
    where $\chi \in \Z$ is a lift of $-u^{-1}w$.
    In particular, $\gamma_u$ is upper-triangular. This then gives $M_\ell = \varepsilon_2D(\ell,1)T^{\chi}$.
\end{remark}

\subsection{Symmetry}\label{subsec: symmetry}

We now focus on conditions for $\Phi_\ell^f$ to be symmetric, leading to \cref{thm: summary conditions} \ref{item: intro symmetry}

\begin{lemma}\label{lem: symmetric poly}
    Let $F\in \C[x,y]$ be an irreducible polynomial of degree greater than~$1$. Then $F$ is a symmetric polynomial if and only if  we have the equality of vanishing loci $V(F(x,y))=V(F(y,x))$.
\end{lemma}
\begin{proof}
    Let $F'(x,y):=F(y,x)$. By the Nullstellensatz and the hypothesis on the loci, we have $\sqrt{\langle F \rangle }= \sqrt{\langle F' \rangle }$.
    Since $F$ and $F'$ are irreducible, this simplifies as
    $\langle F \rangle = \langle F'\rangle$. This means that there exists $\lambda \in \C^\times$ such that $F'(x,y)=\lambda F(x,y)$. Swapping $x$ and $y$, we find $F=\lambda^2F$, meaning that $\lambda=\pm 1$. If $\lambda = -1$, then $F(x,x)=0$. In particular, $F$ is divisible by $x-y$, which contradicts that $F$ is irreducible, since $\deg F>1$. Therefore, $\lambda =1$ and $F$ is symmetric. 
\end{proof}
\begin{theorem}[ = \cref{thm: summary conditions} \ref{item: intro symmetry}]
    Assume that $\Phi_\ell^f$ is denominator-free. Then it is symmetric if and only there exists $\varepsilon \in \{\pm 1\}$ such that $\varepsilon u I \in H$.
\end{theorem}
\begin{proof}
    By \cref{lem: symmetric poly}, $\Phi_\ell^f$ is symmetric if and only if the locus $V(\Phi_\ell^f)$ is symmetric.
    By \cref{thm:properties-Phi-ell-f} \ref{item: classification Phi}, this holds if and only if for any pairs $(E,\alpha),(E',\alpha')\in Y_H^{(-1)}(\C)$, there exists an $\ell$-isogeny $(E,\alpha)\to (E',\alpha')$ exactly when there exists an $\ell$-isogeny $(E',\alpha')\to (E,\alpha)$. For generic $E,E'$, there is no complex multiplication, so, up to $\pm I$, there is at most one $\ell$-isogeny $E'\to E$. In particular, given an $\ell$-isogeny $\phi: E\to E'$, any isogeny $E'\to E$ must be $\pm \hat \phi$. In turn, if $\phi$ respects $\alpha,\alpha'$, then, for any $\varepsilon \in \{\pm 1\}$, we have
    $$
    \alpha \circ \varepsilon\hat \phi = (\varepsilon I)^t\circ \alpha \circ \hat \phi \circ \phi \circ \phi^{-1} = (\varepsilon I)^t\circ\alpha \circ [\ell] \circ \phi^{-1} = (\varepsilon uI)^t\circ \alpha \circ \phi^{-1} = (\varepsilon uI)^t \circ \alpha'.
    $$
    Therefore, $\varepsilon \hat \phi$ respects $\alpha',\alpha$ if and only if $\varepsilon uI\in H$.
\end{proof}

\subsection{Explicit representatives and integrality}\label{subsec: explicit representatives}

In this final subsection, we construct explicit upper-triangular representatives of the cosets $\gamma\in \Gamma_H \backslash \Gamma_H M_\ell \Gamma_H$, under the assumption that $\Phi_\ell^f$ is denominator-free. Then the term $q(\gamma \tau)$ can be expressed in terms of~$q(\tau)$, and we can compute the $q$-expansion of each factor~$f(\gamma \tau)$ in \cref{Phi=Phi tau} from the $q$-expansion of $f$. We also use these formulas to study when $\Phi_\ell^f$ has integer coefficients. The first step is to show that the representatives can stem from the normal subgroup $\Gamma(n) \unlhd \Gamma_H$.

\begin{lemma}\label{lem: principal representative}
    Every coset in $\Gamma_H \backslash \Gamma_H M_\ell \Gamma_H$ admits a representative of the form 
    $$
    D(\ell,1) \gamma \gamma_u, \qquad \text{with }\gamma\in \Gamma(n).
    $$
\end{lemma}
\begin{proof}
    Consider a left coset represented by some $M_\ell \gamma_0$ with $\gamma_0 \in \Gamma_H$. Let $\gamma':= \gamma_u\gamma_0\gamma_u^{-1}$, so that $M_\ell \gamma_0 = D(\ell,1)\gamma' \gamma_u$. We then seek $\gamma \in \Gamma(n)$ and $V\in \Gamma_H$ such that 
    \begin{equation}\label{goal gamma V}
    VD(\ell,1)\gamma' = D(\ell,1)\gamma.
    \end{equation}
    Let $V\in \SL_2(\Z)$ satisfy 
    $$
    \pi_n(V) = D(u,1)\pi_n(\gamma')^{-1}D(u^{-1},1)=M_u^{-1} \pi_n(\gamma_0)^{-1}M_u\qquad \text{and}\qquad \pi_\ell(V) = I.
    $$
    The first condition imposes that $V\in \Gamma_H$. Let $\gamma:=D(1/\ell,1)VD(\ell,1)\gamma'$, so that \cref{goal gamma V} holds. By the second condition on $V$, we have $\ell \mid V_{1,2}$, so $\gamma$ has integer coefficients. Moreover, we have
    $$\pi_n(\gamma) = D(u^{-1},1)\pi_n(V)D(u,1)\pi_n(\gamma')=I,$$
    meaning $\gamma \in \Gamma(n)$.
\end{proof}

\begin{theorem}\label{thm: classification cosets}
    Assume that $\Phi_\ell^f$ is denominator-free; let~$v,w,\varepsilon_1,\varepsilon_2$ be as in \cref{equivalence}. Let $\chi \in \Z$ be a lift of $-u^{-1}w$, and choose $\gamma_u = \varepsilon_2 T^{\chi}$ as in \cref{cor: setting denominator free}. Consider a coset in $\Gamma_H \backslash \Gamma_H M_\ell \Gamma_H$ represented by $D(\ell,1)\gamma \gamma_u$, with $\gamma\in \Gamma(n)$. Then we have:
    \begin{enumerate}
        \item If $\ell \mid \gamma_{2,1}$, then
        $$
        \Gamma_HD(\ell,1)\gamma \gamma_u = \Gamma_H\varepsilon_2 D(\ell,1) T^{\chi}.
        $$
        \item If $\ell \nmid \gamma_{2,1}$, let $b\in \Z$ be such that 
        $$
        \pi_n(b) = v \quad\text{and}\quad \pi_\ell(b) = \pi_\ell(\gamma_{2,1})^{-1} \pi_\ell(\gamma_{2,2}).
        $$
        Then we have
        $$
        \Gamma_H D(\ell,1) \gamma \gamma_u = \Gamma_H \varepsilon_1 \varepsilon_2 D(1,\ell)T^{b+\chi}.
        $$
    \end{enumerate}
    In particular, a set of upper-triangular representatives of the cosets in $\Gamma_H \backslash \Gamma_H M_\ell \Gamma_H $ is given by
    \begin{equation}
    \left\{\varepsilon_2 D(\ell,1)T^{\chi}\right\} \cup \left\{\varepsilon_1 \varepsilon_2 D(1,\ell)T^{b_\nu + \chi}: \nu \in \Z/\ell \Z\right\},
    \label{eq: list of representatives}
    \end{equation}
    where $b_\nu\in \Z$ is any number satisfying $b_\nu \equiv v \mod n$ and $b_\nu \equiv \nu \mod \ell$.
\end{theorem}
\begin{proof}
    Assume that $\ell \mid \gamma_{2,1}$. Let $V = D(\ell,1)\gamma D(1/\ell,1)$, which then has integer entries. The matrix $V$ also satisfies $\det V = 1$, so $V \in \SL_2(\Z)$. Moreover, we also have $\pi_n(V) = I$, so $V\in \Gamma(n) \subseteq \Gamma_H$. In particular, from $D(\ell,1)\gamma = VD(\ell,1)$, substituting $\gamma_u$ with the expression in \cref{cor: setting denominator free}, we obtain the result.
    Assume now that $\ell \nmid \gamma_{2,1}$ and let $b$ be as in the statement. Define
    $$V = \varepsilon_1 D(\ell,1)\gamma T^{-b}D(1,1/\ell).$$
    We have
    $$
    \pi_\ell((\gamma T^{- b})_{2,2}) = -\pi_\ell(\gamma_{2,1})\pi_\ell(b) +\pi_\ell(\gamma_{2,2})= 0.
    $$
    Therefore, $V$ has integer coefficients. Moreover, we also have $\det V = 1$, so $V\in \SL_2(\Z)$. In addition, we compute
    $$
    \pi_n(V) = \varepsilon_1 D(u,1)\pi_n(T)^{-\pi_n(b)}D(1,u^{-1}) = \varepsilon_1 U(u,v)^{-1}\in H.
    $$
    Thus, $V\in \Gamma_H$. Ultimately, by construction we have
    $$
    D(\ell,1)\gamma = V \varepsilon_1 D(1,\ell) T^{b},
    $$
    which concludes that $\varepsilon_1\varepsilon_2D(1,\ell)T^{b+\chi}$ also represents the class of $D(\ell,1)\gamma \gamma_u$.
    By \cref{lem: principal representative}, every coset admits a representative of the form $D(\ell,1)\gamma \gamma_u$ with $\gamma \in \Gamma(n)$. As a consequence, every coset is represented by an element in \cref{eq: list of representatives}. Since we have $|\Gamma_H \backslash \Gamma_H M_\ell \Gamma_H| = \ell +1$ and the cardinality of \cref{eq: list of representatives} is at most $\ell+1$, this concludes that both sets have cardinality $\ell+1$, which finishes the proof.
\end{proof}
\begin{corollary}\label{cor: expression product}
    Let $n_0\mid n$ be the width of $\Gamma_H$, and $q_{n_0} = \exp\left(2\pi i \tau\right/n_0)$.
    Let $\tilde f \in \C((q_{n_0}))$ be the $q$-expansion of $f$ at infinity, so that $f(\tau)=\tilde f(q_{n_0})$.
    Then, continuing with the notation from \cref{thm: classification cosets}, the $q$-expansions of the functions $f(\gamma \tau)$ with $\gamma \in \Gamma_H \backslash \Gamma_H M_\ell \Gamma_H$ at infinity are given by the evaluation of $\tilde f$ at the variables
   \begin{equation}
       q_{n_0}(\ell(\tau+\chi))=\zeta_{n_0}^{-w} q_{n_0}^{\ell}, \qquad 
       q_{n_0}\left(\frac{1}{\ell}(\tau+b_{n_0\nu-\chi}+\chi)\right)=\zeta_{n_0}^{u^{-1}(v-u^{-1}w)}\zeta_\ell^\nu q_{n_0}^{1/\ell},
    \label{eq: q expansion variables}
    \end{equation}
    where $\nu \in \Z/\ell \Z$.
\end{corollary}
\begin{proof}
    It suffices to use the set of representatives from \cref{thm: classification cosets}, and to compute their action on $\tau$. 
\end{proof}
\begin{corollary}\label{cor: coefficients in Z}
    Assume that
    \begin{enumerate}
        \item $\Phi_\ell^f$ is denominator-free,
        \item one can take $v=w=0$ in \cref{equivalence},
        \item the $q$-expansion of $f$ has integer coefficients and leading coefficient 1.
    \end{enumerate}
    Then $\Phi_\ell^f(x,y) \in \Z[x,y]$.
\end{corollary}
\begin{proof}
    For $0\leq r \leq \ell+1$, let $S_{\ell+1-r}$ denote the $(\ell+1-r)$-th elementary symmetric polynomial.
    By \cref{thm:properties-Phi-ell-f} \ref{item: Phi=Phi tau}, we have
    $$
    S_{\ell+1-r}(f(\gamma\tau)_{\gamma \in \Gamma_H \backslash \Gamma_H M_\ell \Gamma_H}) = R_r(f(\tau)),
    $$
    where $R_r(x)\in K_H[x]$ is the $r$-th coefficient of $\Phi_\ell^f$ in $y$.
    Using \cref{cor: expression product} to write the $q$-expansions of both sides, since $v=w=0$, $\zeta_{n_0}$ does not appear in the left-hand side. In addition, by Galois invariance, $\zeta_\ell$ does not appear on the left-hand side either. Therefore, $R_r(x)\in \Q[x]$. Moreover, since the $q$-expansion of $f$ has integer coefficients, so does the $q$-expansion of the left-hand side. As the leading coefficients of the $q$-expansion of $f$ is 1, we conclude by Hasse's $q$-expansion principle (see \cite[Exercise 11.12]{Cox1997}, replacing $j$ with $f$).
\end{proof}
\begin{remark}\label{rem: computation from q exp}
    This method consisting of comparing the $q$-expansions of both sides provides a quadratic-time algorithm to compute $R_r$ and hence $\Phi_\ell^f$; it also provides an additional proof that $\Phi_\ell^f$ has coefficients in $\Q(\zeta_n)$.
\end{remark}

\section{Heights of modular polynomials}\label{section heights}

In this section, we prove that modular polynomials arising from genus~$0$ covers of~$X(1)$ have smaller height, as stated in \cref{thm:main-height}. Our proof uses an evaluation-interpolation approach, using similar tools as \cite{kieffer:hal-03226568}.

\subsection{Reminders on heights}
\label{subsec:heights}

The absolute logarithmic height is a real-valued function which quantifies the ``size'' of various objects obtained from a number field, such as varieties over it or polynomials with coefficients in it.

Let $K$ be a number field. Following \cite[§B.1]{hindry2000diophantine}, let $M_K$ be the set of absolute values on $K$ whose restriction to $\Q$ is either $|\cdot|_\infty$ or $p$-adic, $M_K^\infty\subseteq M_K$ be the set of those that are archimedean and $M_K^0\subseteq M_K$ the set of those that are nonarchimedean. Consider a polynomial $F\in K[x_1,\ldots,x_l]$, and write  $F(x_1,\ldots,x_l)=\sum_{\omega\in \Omega}a_\omega x_1^{\omega_1}\ldots x_l^{\omega_l}$, so that we view $F$ as the affine vector $(a_\omega)_{\omega\in \Omega}$. Let $v\in M_K$, and define
$$
|F|_v := |(a_\omega)_{\omega\in \Omega}|_v:=\max_{\omega\in \Omega}|a_\omega|_v.
$$
\begin{lemma}\label{lem:a v polynomials}
    The nonarchimedean and archimedean inequalities extend to polynomials. More precisely, given a finite set of polynomials $\mathcal X\subseteq  K[x_1,\ldots,x_l]$, we have
    $$
    \left|\sum_{F\in \mathcal X} F\right|_v \leq \begin{cases}
        \max_{F\in \mathcal X}|F|_v \qquad &\text{if $v$ is nonarchimedean} \\
        \sum_{F\in \mathcal X} |F|_v   &\text{if $v$ is archimedean}.
    \end{cases}
    $$
\end{lemma}
\begin{proof}
    This is a direct calculation from the definition.
\end{proof}
For each $v\in M_K$, let $d_v:=[K_v:\Q_v]$, where the subscripts denote the completion with respect to $v$. Also, define $d_K := [K:\Q]$.
Following \cite[§B.7]{hindry2000diophantine}, we set
$$
h(F) := h[1:\ldots:a_\omega:\ldots ]_{\omega\in \Omega}:=\sum_{v\in M_K} \frac{d_v}{d_K} \log \max(1,|F|_v).
$$
In particular, when $F$ has coefficients in $\Z$ and is nonzero, we recover the intuitive expression
$$
h(F) = \log \max_{\omega\in \Omega} |a_\omega|.
$$

\subsection{Height bounds from evaluations}
Choose an integer $N$ and a set of distinct complex numbers $\{\tau_1,\ldots,\tau_N\}\subseteq \UH$ whose image under $f$ is $\mu_N(\C)$, the set of $N$-th roots of unity, which we enumerate as $\{\xi_1,\ldots, \xi_N\}$.

\begin{lemma}\label{lem: delta}
    Let $1\leq m\leq N$, consider a set of indices $I\subseteq \{1,\ldots, N\}$ such that $|I|=m$ and define for any $\iota \in  I$:
    $$\delta_\iota=\prod_{s \in I\setminus \{\iota\}}(\xi_{\iota}-\xi_s).$$
    For any $v\in M_{\Q(\mu_N(\C))}$, we have the bounds
    $$
    |\delta_\iota|_v \geq \begin{cases}
        |N|_v & \text{ if } v \text{ is nonarchimedean} \\
        N\cdot2^{-(N-m)}  & \text{ if } v \text{ is archimedean}.
    \end{cases}
    $$
\end{lemma}
\begin{proof}
We have
\begin{equation}\label{delta}
\delta_\iota = \prod_{\substack{1\leq s \leq N \\ s \neq \iota}}(\xi_{\iota}-\xi_s)\prod_{\substack{1\leq s \leq N \\ s \notin I}}(\xi_{\iota}-\xi_s)^{-1}.
\end{equation}
Since $t^N-1 = \prod_{1\leq s \leq N}(t-\xi_s)$, the first product is equal to the derivative $(t^N-1)'=Nt^{N-1}$ evaluated at $t=\xi_\iota$. Consequently, the norm of the first product is equal to $|N\xi_{\iota}^{N-1}|_v= |N|_v$.
If $v$ is nonarchimedean, then each factor $|\xi_{\iota}-\xi_s|_v$ in the second product in \cref{delta} is smaller than or equal to $\max(|\xi_\iota|_v, |\xi_s|_v)=1$. This implies that
$$
|\delta_\iota|_v \geq |N|_v.
$$
On the other hand, if $v$ is archimedean, using \cref{delta} again, we simply have the bound $|\xi_{\iota}-\xi_s|_v \leq 2$, which yields
\begin{equation}
|\delta_\iota|_v \geq |N|_v\prod_{\substack{1\leq s \leq N \\ s \notin I}}\frac{1}{2} = N \frac{1}{2^{N-m}}. \tag*{\qedhere}
\end{equation}
\end{proof}

\begin{proposition}\label{prop1}
    Let $K$ be a number field, let $F(x,y)\in K[x,y]$ be a polynomial, let $N\geq \deg_x (F)+1$, and write $L=K(\mu_N(\C))$. Fix an integer $m$ satisfying $\deg_x(F)+1\leq m\leq N$, and let $I\subseteq \{1,\ldots,N\}$ be a set of $m$ distinct indices. For any $1\leq \iota  \leq N$, let $f_\iota(y) := F(\xi_\iota,y)$. Then, for any absolute value $v\in M_{L}$, we have
    $$
    |F|_v \leq C_v \max_{\iota\in I} |f_\iota|_v,
    $$
    where 
    $$C_v=\begin{cases}1/|N|_v& \text{if $v$ is nonarchimedean}\\
    \frac{m}{N}2^{N-m} \begin{pmatrix}
    m-1 \\ \lceil(m-1)/2\rceil
    \end{pmatrix}&\text{if $v$ is archimedean.}\end{cases}$$
    In particular, we have $C_v \geq 1$.
\end{proposition}
\begin{proof}
We interpolate $F(x,y)$ with respect to $x$ using $m$ values of $f_\iota(y)$. We have
$$
F(x,y) = \sum_{\iota\in I} f_{\iota}(y) W_\iota(x), \qquad
\text{ where }\qquad W_\iota(x) = \prod_{s \in I\setminus \{\iota\}} \frac{x-\xi_s}{\xi_\iota-\xi_s}.$$
Let $\delta_\iota = \prod_{s \in I\setminus\{\iota\}}(\xi_\iota-\xi_s)$ and let $\theta_{\iota,r}\in L$ be the coefficient of $x^r$ in the numerator of $W_\iota(x)$. We then have
$$
W_\iota(x) = \frac{1}{\delta_\iota} \sum_{r=0}^{m-1}\theta_{\iota,r} x^r \qquad\text{ and }\qquad F(x,y) = \sum_{\iota\in I} \sum_{r=0}^{m-1} f_\iota(y) \frac{1}{\delta_\iota} \theta_{\iota,r}x^r.
$$
Therefore, for any $v\in M_L$,
\begin{equation}\label{v norm phi}
|F|_v = \max_{0\leq r \leq m-1} \left| \sum_{\iota\in I} f_{\iota}(y) \frac{1}{\delta_\iota} \theta_{\iota,r} \right|_v.
\end{equation}
In addition, we can expand the expression of $W_\iota$ in order to obtain a value for $\theta_{\iota,r}$, namely
\begin{equation}\label{alpha}
\theta_{\iota,r} =(-1)^{m-1-r}\sum_{\substack{J\subseteq I\setminus  \{\iota\} \\|J|=m-1-r}}  \prod_{t\in J}\xi_{t}.
\end{equation}
We now distinguish the two classical cases for $v$.

\subsubsection*{Nonarchimedean case} Using the nonarchimedean inequality from \cref{lem:a v polynomials}, for any $0\leq r\leq m-1$, the expression in the maximum in \cref{v norm phi}
is smaller than or equal to 
$$\max_{\iota\in I} \left| f_\iota(y) \frac{1}{\delta_\iota} \theta_{\iota,r} \right|_v.$$
Using the nonarchimedean inequality again, we get from \cref{alpha}
$$
|\theta_{\iota,r}|_v \leq \max_{\substack{J\subseteq I\setminus  \{\iota\} \\|J|=m-1-r}} \left| \prod_{t\in J}\xi_t \right|_v =1.
$$
Using \cref{lem: delta} in the nonarchimedean case, we then obtain
$$
|F|_v \leq \max_{0 \leq r\leq m-1} \max_{\iota\in I} |f_\iota(y)|_v \frac{1}{|\delta_\iota|_v} |\theta_{\iota,r}|_v  \leq  \frac{1}{|N|_v} \max_{\iota\in I}|f_{\iota}(y)|_v,
$$
which is the result we were seeking in the nonarchimedean case. By definition of $|N|_v$, we directly find $C_v\geq 1$ in this case.

\subsubsection*{Archimedean case}
Getting back to \cref{v norm phi}, using the archimedean inequality from \cref{lem:a v polynomials}, we can now bound the right-hand side by
$$
\max_{0\leq r \leq m-1} \sum_{\iota\in I} |f_{\iota}(y)|_v \frac{1}{|\delta_\iota|_v} |\theta_{\iota,r}|_v \leq 
m \max_{0\leq r \leq m-1} \max_{\iota\in I} |f_{\iota}(y)|_v \frac{1}{|\delta_\iota|_v} |\theta_{\iota,r}|_v.
$$
For $\theta_{\iota,r}$, again from \cref{alpha}, we have 
$$|\theta_{\iota,r}|_v \leq \sum_{\substack{J\subseteq I\setminus  \{\iota\} \\|J|=m-1-r}} \left| \prod_{t\in J}\xi_{t} \right|_v = \begin{pmatrix}
    m-1 \\ r
\end{pmatrix}.$$
Using \cref{lem: delta} in the archimedean case, we get
\begin{align*}
|F|_v &\leq m \max_{0\leq r \leq m-1} \max_{\iota\in I} |f_{\iota}(y)|_v \frac{1}{N} 2^{N-m} \begin{pmatrix}
    m-1 \\ r
\end{pmatrix}\\
&= \frac{m}{N}2^{N-m} \begin{pmatrix}
    m-1 \\ \lceil(m-1)/2\rceil
\end{pmatrix}  \max_{\iota\in I} |f_{\iota}(y)|_v.
\end{align*}
The binomial coefficient is larger than or equal to 1 by definition. On the other hand, taking $t=N-m$, the binomial expansion gives
$$
2^t = (1+1)^t \geq 1 +t \geq 1+t/m = N/m,
$$
proving that $C_v\geq 1$ also in this case.
\end{proof}

\begin{corollary}\label{prop}
    Continuing with the notation from \cref{prop1}, we have the following bound on the height of $F$ :
    \begin{equation}\label{h(F)}
    h(F) \leq \frac{N}{N-m+1} \max_{1 \leq \iota  \leq N} h(f_\iota) +\frac{1}{2} \log\left(\frac{m^2}{2m-1} \right) +N\log 2-\frac{1}{2}\log \pi.
    \end{equation}
\end{corollary}
\begin{proof}
    Fix some $v\in M_L$ and consider the result of \cref{prop1}. Since every subset $I\subseteq \{1,\ldots,N\}$ of $m$ indices contains at least one $\iota \in  I$ with $|F|_v\leq C_v |f_\iota|_v$, there can be at most $m-1$ indices $\iota \in  \{1,\ldots,N\}$ with $|F|_v> C_v |f_\iota|_v$. Hence, there must be at least $N-m+1$ distinct indices $\mathsf{j}_1,\ldots, \mathsf{j}_{N-m+1}\in \{1,\ldots,N\}$, possibly dependent on $v$, satisfying 
$$
|F|_v \leq C_v |f_{\mathsf{j}_t}|_v
$$
for all $1\leq t \leq N-m+1$.
From \cref{prop1}, since $C_v\geq 1$, we deduce that
$$
\max(1,|F|_v) \leq C_v \max(1,|f_{\mathsf{j}_t}|_v).
$$
Multiplying these inequalities for all $1\leq t\leq N-m+1$, we obtain
$$
\max(1,|F|_v)^{N-m+1} \leq C_v^{N-m+1} \prod_{t=1}^{N-m+1}\max(1,|f_{\mathsf{j}_t}|_v).
$$
Taking logarithms yields
\begin{align*}
    (N-m+1) \log \max(1,|F|_v) &\leq (N-m+1) \log C_v + \sum_{t=1}^{N-m+1}\log \max(1,|f_{\mathsf{j}_t}|_v) \\
    &\leq (N-m+1) \log C_v + \sum_{\iota=1}^{N}\log \max(1,|f_{\iota}|_v).
\end{align*}
Averaging over all standard absolute values then gives us
\begin{align*}
    h(F) &\leq  \sum_{v\in M_{L}} \frac{d_v}{d_L} \log C_v+ \frac{1}{N-m+1} \sum_{v\in M_L} \frac{d_v}{d_L} \sum_{\iota=1}^{N}\log \max(1,|f_{\iota}|_v)\\
    &=C+ \frac{1}{N-m+1}\sum_{\iota=1}^{N}h(f_{\iota}) \leq C+\frac{N}{N-m+1} \max_{1\leq \iota  \leq N} h(f_{\iota})
\end{align*}
where $C = \sum_{v\in M_{L}} \frac{d_v}{d_L} \log C_v$. Let us compute $C$. When $v$ is archimedean, $C_v=C_{\infty}$ does not depend on $v$, so 
$$
\sum_{v\in M_L^\infty}\frac{d_v}{d_L}\log C_v=\log C_{\infty}\sum_{v\in M_L^\infty}\frac{d_v}{d_L}=\log C_{\infty}.
$$
On the other hand, when $v$ is nonarchimedean, one has
\begin{align*}
    \sum_{v\in M_L^0}\frac{d_v}{d_L}\log C_v &= \sum_{v\in M_L^0}\frac{d_v}{d_L}\log \frac{1}{|N|_v} = \sum_{p\text{ prime}}\sum_{\substack{v\in M_L^0 \\ v|p}}\frac{d_v}{d_L}\ord_p(N) \log p \\
    &=\sum_{p|N}\ord_p(N) \log p\sum_{\substack{v\in M_L^0 \\ v|p}}\frac{d_v}{d_L}=\log N.
\end{align*}
Therefore, $$C=\log N +\log C_\infty=\log m +(N-m)\log 2+\log\begin{pmatrix}
    m-1 \\ \lceil(m-1)/2\rceil
    \end{pmatrix}.$$
    Finally, using the Stirling formula from \cite{Junod2022CoefficientBinomialCentral}, we have $$\log\begin{pmatrix}
    m-1 \\ \lceil(m-1)/2\rceil
    \end{pmatrix}  \leq  m\log 2 -\frac{1}{2} \log(2m-1)-\frac{1}{2}\log \pi.$$
Substituting,
\begin{equation}
 C\leq\log m-\frac{1}{2}\log(2m-1)+N\log 2-\frac{1}{2}\log \pi. \tag*{\qedhere}
\end{equation}
\end{proof}

\subsection{Height bounds for modular polynomials}\label{subsec: main thrm}
Applying \cref{prop} to modular polynomials allows us to prove \cref{thm:main-height}. Let~$X, K, f, p, d, \ell$, and $\Phi_\ell^f$ be as in the statement of the theorem (which is also repeated as \cref{thm+proof: heights} below).
\begin{lemma}\label{lem: height covering}
    There exists a constant $c>0$ such that for any $P \in X(\overline \Q)$:
$$
\left|dh(f(Q))-h(j\circ p(Q))\right| \leq c.
$$
\end{lemma}
\begin{proof}
    It suffices to apply \cite[B.2.5. (b)]{hindry2000diophantine}, to the degree-$d$ morphism $j\circ p\circ f^{-1}:\Proj^1_{K}\to \Proj^1_{K}$ at the point $f(Q)$.
\end{proof}

\begin{lemma}\label{lem: height isogeny}
    Let $Q,Q'\in X(\overline \Q)$, be such that $(p(Q),p(Q'))$ is either $(\infty,\infty)$ or corresponds to a pair of $\ell$-isogenous elliptic curves. Then we have
    $$
    dh(f(Q)) \leq dh(f(Q')) + 2c + \Upsilon_1(h(f(Q'))+c,\ell),
    $$
    where $\Upsilon_1(a,\ell):=2\sqrt3\sqrt{a+6\log(6 \ell)-5}+6\log (6\ell)$.
\end{lemma}
\begin{proof}
    Let $P=p(Q)$ and $P'=p(Q')$.    
    If $(P,P')$ corresponds to a pair of $\ell$-isogenous elliptic curves, by \cite[Theorem 1.1 (2)]{doi:10.1142/S1793042119500295}, we have
\begin{equation}
h(j(P))-6\log(1+h(j(P)) \leq h(j(P'))+ 6\log \ell.
\label{eq: Pazuki}
\end{equation}
On the other hand, if $(P,P') = (\infty,\infty)$ then $j(P)=j(P')=\infty$, so \cref{eq: Pazuki} also holds in this case.
This inequality is equivalent to
\begin{equation}\label{W-1}
h(j(P)) \leq -6W_{-1}\left(-\frac{1}{6\ell}e^{-(h(j(P'))+1)/6} \right)-1,
\end{equation}
where $W_{-1}$ denotes the lower branch of the Lambert function. In particular, it satisfies for any $0<x <e^{-1}$:
\begin{equation}\label{Cha}
W_{-1} (-x)>\log x-\sqrt{-2(\log x+1)},
\end{equation}
as shown in \cite[theorem 1]{Chatzigeorgiou}. For any $\ell \geq 1 > \frac{1}{6}e^{5/6}$ and $a\geq 0$, we have $0<\frac{1}{6\ell}e^{-(a+1)/6}<e^{-1}$. We can then apply \cref{Cha} to the right hand side of \cref{W-1}, and we obtain
$$
-6W_{-1}\left(-\frac{1}{6\ell}e^{-(a+1)/6} \right)-1 <a+ \Upsilon_1(a,\ell),
$$
where $\Upsilon_1(a,\ell)$ is defined as above.
Therefore, using \cref{lem: height covering}, we have
\begin{align*}
dh(f(Q))-c\leq h(j(P)) &\leq h(j(P'))+\Upsilon_1(h(j(P')),\ell)\\
&\leq dh(f(Q'))+ c+\Upsilon_1(h(f(Q'))+c,\ell). \qedhere
\end{align*}
\end{proof}

\begin{theorem}[= \cref{thm:main-height}]\label{thm+proof: heights}
Let $X$ be a smooth, projective curve of genus zero defined over a number field~$K$ endowed with a covering map $p: X\to X(1)$ of degree~$d$, and choose an isomorphism 
    $f:X \overset{\sim}{\longrightarrow} \mathbb{P}^1_{K}$. Let~$\ell$ be a prime. Suppose that there exists a polynomial $\Phi_\ell^f \in \overline \Q[x,y]$ such that:
    \begin{enumerate} 
        \item
        $\Phi_\ell^f$ has degree $\ell+1$ in both variables,
        and is monic as a polynomial in $y$.
        \item 
        For every pair $(Q,Q')\in 
        X(\overline \Q) \times X(\overline \Q)$ 
        satisfying
        \[
            \Phi_\ell^f(f(Q),f(Q')) = 0,
        \]
        the pair $(p(Q),p(Q'))$ is either $(\infty,\infty)$ or corresponds to~$\ell$-isogenous elliptic curves.
    \end{enumerate}
    Then we have
    \[
    h(\Phi_\ell^f) \leq \frac{6}{d}\ell\log\ell + O_{p,f}\left(\ell\sqrt{\log(\ell)}\right).
    \]
    The implied constants in this estimate are independent of~$\ell$ and depend only on~$p,f$.
\end{theorem}
\begin{proof}
We first apply \cref{prop} to $\Phi_\ell^f$. Let $Q_1,\ldots, Q_N \in X(\overline \Q)$ be a set of distinct points whose image under $f$ us $\mu_N(\C)$. The bound from \cref{h(F)} increases with the number $m$ of interpolation points. We therefore take $m$ as small as possible, i.e. $m=\ell+2$.
We then get
$$
h(\Phi) \leq\frac{N}{N-\ell-1} \max_{1 \leq \iota \leq N}h(\Phi(f(Q_\iota),y))+N\log 2+\Upsilon_2(\ell),$$
where $\Upsilon_2(\ell)=\frac{1}{2} \log\left((\ell+2)^2/(2\ell+3) \right) -\frac{1}{2}\log \pi$.

For any $1\leq \iota \leq N$, consider the polynomial $\Phi_\ell^f(f(Q_\iota),y)$. By assumption, it is monic and has coefficients in $\overline \Q$. Its zeroes are $\ell+1$ algebraic numbers that can be written as $f(Q)$ for $Q\in X(\overline \Q)$.
Let $\mathcal{T}_\iota$ be a multiset of such elements (meaning it takes into account the potential multiplicities), so that
$$
\Phi(f(Q_\iota),y)=\prod_{Q \in \mathcal{T}_\iota}(y-f(Q)).
$$
Using \cite[Proposition B.7.2.]{hindry2000diophantine}, we have
$$
    h(\Phi(f(Q_\iota),y)) = h\left(\prod_{Q\in \mathcal{T}_\iota}(y-f(Q)) \right) \leq \sum_{Q \in \mathcal{T}_\iota}h(f(Q))+2\log 2\cdot(\ell+1).
$$
Since $h(f(Q_\iota)) = 0$,  \cref{lem: height isogeny} gives for each $Q \in \mathcal{T}_\iota$:
$$
dh(f(Q)) \leq 2c + \Upsilon_1(c,\ell).
$$
In addition, we have $|\mathcal{T}_\iota|=\deg_x\Phi=\ell+1$.
Thus, 
$$
h(\Phi_\ell^f(f(Q_\iota),y)) \leq (\ell+1)\left(\frac{1}{d}(2c+\Upsilon_1(c,\ell))+2\log 2\right)=:\Upsilon_3(c,\ell).
$$
This implies that
\begin{equation}\label{Upsilon_4}
h(\Phi_\ell^f) \leq \frac{N}{N-\ell-1}\Upsilon_3(c,\ell)+N\log 2+\Upsilon_2(\ell)=:\Upsilon_4(N,\ell).
\end{equation}
We now want to choose $N$ in order to minimize $\Upsilon_4(N,\ell)$.
Taking the second derivative with respect to $N$, we see that $\Upsilon_4(N,\ell)$ is convex. Taking the first derivative, we see that it has a global minimum at
$$
N= \ell+1+ \Upsilon_5(c,\ell),
$$
where $\Upsilon_5(c,\ell):=\sqrt{(\ell+1)\Upsilon_3(c,\ell)/\log(2)}$. However, the inequality \cref{Upsilon_4} only holds for $N$ integral.
We then have
$$
h(\Phi_\ell^f) \leq \Upsilon_4(\ell+1+\lceil \Upsilon_5(c,\ell)\rceil,\ell) \leq \Upsilon_4(\ell+2+\Upsilon_5(c,\ell),\ell).
$$
We then compute 
\begin{align*}
\Upsilon_4(\ell+2 +\Upsilon_5(c,\ell), &\ell) = \frac{6}{d}\ell\log\ell\\
&+\frac{2\sqrt3}{d}(\ell+1)
\sqrt{c+6\log(6\ell)-5}
+2(\log2)\,\Upsilon_5(c,\ell)\\
&+(\ell+1)\left(
\frac{2c+6\log6}{d}+3\log2
\right)\\
&+\left(\frac6d+\frac12\right)\log\ell
-\frac12\log(2\pi)\\
&+\log\left(1+\frac2\ell\right)
-\frac12\log\left(1+\frac{3}{2\ell}\right)
+\frac{\log2}{1+\Upsilon_5(c,\ell)}.
\end{align*}
Notice that $\Upsilon_1(c,\ell)=O(\log \ell)$, $\Upsilon_3(c,\ell)=O(\ell \log \ell)$ and $\Upsilon_5(c,\ell)=O(\ell \sqrt{\log \ell})$, where the implied constants are allowed to depend on~$p$ and~$f$. Thus the second term is $O(\ell\sqrt{\log(\ell)})$.
\end{proof}

In fact, we can make the implied constant in \cref{thm:main-height} more precise.

\begin{corollary}
\label{cor:explicit}
    Let~$c$ be as in \cref{lem: height covering}. Then in \cref{thm+proof: heights}, we have
    $$
    h(\Phi_\ell^f) \leq 
\frac{6}{d}\ell\log\ell
+60\left(1+\frac{1}{d}\right)(1+c)
\ell\sqrt{\log\ell}.
    $$
\end{corollary}
\begin{proof}
For $\ell\geq 2$ and $c\geq 0$, we have
$$
\sqrt{c+6\log(6\ell)-5}
\leq 4(1+c)\sqrt{\log\ell}.
$$
Consequently,
$$
2c+2\sqrt{3}\sqrt{c+6\log(6\ell)-5}
+6\log 6
\leq 30(1+c)\sqrt{\log\ell}.
$$
It follows that
$$
\frac{2c+\Upsilon_1(c,\ell)}{d}+2\log 2
\leq
\left(\frac{43}{d}+2\right)
(1+c)\log\ell,
$$
and hence
$$
2(\log 2)\,\Upsilon_5(c,\ell)
\leq
\left(\frac{11}{\sqrt d}+3\right)
(1+c)(\ell+1)\sqrt{\log\ell}.
$$
Using the exact expression for
$\Upsilon_4(\ell+2+\Upsilon_5(c,\ell),\ell)$
and the inequalities
$$
\Upsilon_2(\ell)\leq\frac12\log\ell,
\qquad
\ell+1\leq\frac32\ell,
$$
we obtain
\begin{align*}
&\Upsilon_4\bigl(\ell+2+\Upsilon_5(c,\ell),\ell\bigr)
-\frac6d\ell\log\ell\\
&\quad\leq
\left(
\frac{51}{d}
+\frac{33}{2\sqrt d}
+\frac{21}{2}
\right)
(1+c)\ell\sqrt{\log\ell}\\
&\quad\leq
60\left(1+\frac1d\right)
(1+c)\ell\sqrt{\log\ell},
\end{align*}
where the last inequality follows from
$d^{-1/2}\leq(1+d^{-1})/2$.
\end{proof}

\section{Examples}\label{section: examples}

This section aims to present examples of our construction of modular polynomials. These will concern two invariants, namely the $\gamma_2$ and the Weber $\f$ modular functions. 
We begin by explaining how these functions fit in our framework. 

\subsection{Modular functions as morphisms}\label{section: starting thw Gamma}

Let~$H$ be a subgroup of~$\GL_2(\Z/n\Z)$ for some~$n\geq 1$.
The following proposition recalls how a weight~0 modular function $f$ of level $\Gamma_H = \SL_2(\Z)\cap \pi_n^{-1}(H)$ can be considered as an algebraic morphism $X_H^{(-1)}\to \Proj^1$.

\begin{proposition}\label{prop: functions GAGA}
    Let $f$ be a nonconstant modular function of weight 0 for $\Gamma_H$. Then we have:
    \begin{enumerate}
        \item $f$ defines a morphism of algebraic curves $X_{H}^{(-1)}\to \Proj^1$ over $\C$.
        \item \label{item: q expansion principle} For any number field~$K$ containing $K_H$, the morphism $f:X_{H}^{(-1)} \to \Proj^1$ is defined over $K$ if and only if the $q$-expansion of $f$ at infinity has coefficients in $K$.
        \item \label{item: Hauptmodul}
        The morphism $f:X_H^{(-1)}\to\Proj^1$
        is an isomorphism if and only if there exists a rational
        fraction $R\in\C(x)$ such that
        $$
        j\circ p=R\circ f
        \qquad\text{and}\qquad
        \deg R=\deg p,
        $$
        where $p:X_H^{(-1)}\to X(1)$ is the natural covering map.
    \end{enumerate}
\end{proposition}
\begin{proof}
We prove each item in order.
\begin{enumerate}
\item As $f$ is invariant under $\Gamma_H$
and is meromorphic on the upper half-plane and at every
cusp, it induces a holomorphic function 
$$\Gamma_H \backslash \UH^* \simeq X_H^{(-1)}(\C) \to \Proj^1(\C).$$ By \cite[§19, Proposition 15]{GAGA}, this map is the analytification of an algebraic morphism $X_H^{(-1)}\to \Proj^1_\C$.

\item This is a reformulation of the results from \cite[Section V, Theorem 3.9, Corollary 3.13 and §4.6-4.7]{Deligne.Rapoport1973}.

\item Suppose $R$ exists with the desired properties.
Since all morphisms are nonconstant morphisms
between smooth projective curves and $j$ is an isomorphism, the degrees satisfy
$$
\deg p=\deg R\cdot \deg f.
$$
Consequently,
$$
\deg p=\deg R
\iff
\deg f=1 \iff \text{ $f$ is an isomorphism},
$$
where the last equivalence holds because we are working with smooth projective
irreducible curves.
Conversely, if $f$ is an isomorphism, then
\[
R=j\circ p\circ f^{-1}:\mathbb P^1\to\mathbb P^1
\]
is a rational function satisfying $j\circ p=R \circ f$, and
$\deg R=\deg p$. \qedhere
\end{enumerate}
\end{proof}

We now fix a congruence subgroup $\Gamma\leq \SL_2(\Z)$ and a modular function $f$ with stabilizer $\Gamma$ and weight 0. 
To apply our formalism, we will also fix a subgroup $H$ satisfying $\Gamma_H = \Gamma$ and consider $f$ as an invariant on $X_H^{(-1)}$. Note that any such~$H$ satisfies $H\leq N_{\GL_2(\Z/n\Z)}(\pi_n(\Gamma))$.

    Our construction of~$\Phi_\ell^f$ in~§\ref{section: construction} relies on the hypothesis that for each prime degree $\ell$ under consideration, we have 
$$(\ell \mod n) \in \det H.$$
Ideally, the choice of $H$ would not depend on $\ell$. In other words, one would hope that the set
\begin{equation}\label{eq: maximum determinant}
\{\det H: H \leq \GL_2(\Z/n\Z), \Gamma_H = \Gamma\}
\end{equation}
admits a maximum for the inclusion.
However, this is not always true, as shown by the following example, found using Magma \cite{BCP97}.

\begin{example}
Let $n=8$, and define the matrices modulo 8 $$A=\begin{pmatrix}
    1 & 7 \\ 1 & 0
\end{pmatrix}, \hspace{5mm }X= \begin{pmatrix}
    1 & 4 \\ 4 & 5
\end{pmatrix},\qquad Y= \begin{pmatrix}
    0 & 1 \\ 1 & 0
\end{pmatrix}.$$
These satisfy $A\in \SL_2(\Z/8\Z)$ and $X,Y\in \GL_2(\Z/8\Z)$.
Also, define the subgroups $$G=\langle A\rangle, \qquad H_X=\langle A,X \rangle, \qquad H_Y = \langle A, Y \rangle.$$
Then we have
$$
H_X \cap \SL_2(\Z/8\Z) = G \qquad \text{and} \qquad H_Y \cap \SL_2(\Z/8\Z) = G.
$$
However, there exists no subgroup $H\leq \GL_2(\Z/8\Z)$ such that 
$$H\cap \SL_2(\Z/8\Z)=G \qquad \text{and} \qquad\det H_X, \det H_Y \leq \det H.$$
\end{example}

In both our examples below, the set~\eqref{eq: maximum determinant} indeed admits a maximum (in fact, we can choose~$H$ such that $\det H = (\Z/n\Z)^\times$).

\subsection{The \texorpdfstring{$\gamma_2$}{} modular function}

Consider first the \textit{Dedekind eta function} $$\eta(\tau) := q^{1/24} \prod_{n=1}^\infty (1-q^n),$$
where $q = e^{2\pi i \tau}$. 
Let $\gamma_2$ be the function on $\UH$ defined by
$$
\gamma_2(\tau) = \frac{3g_2(\tau)}{4\pi^4 \eta(\tau)^8},\qquad\text{where}\qquad g_2(\tau)=60 \sum_{(m,n)\in \Z^2 \setminus \{(0,0)\}} (m+n\tau)^{-4}.
$$
As shown in \cite[§12.A]{Cox1997}, we have $\gamma_2^3=j$. The largest congruence subgroup $\Gamma_{\gamma_2} \leq \SL_2(\Z)$  leaving $\gamma_2$ invariant can be described as follows. Consider the quaternionic subgroup of $\SL_2(\Z/3\Z)$ given by
\[G_{\gamma_2}=
\left\{\pm\begin{pmatrix}
    1 & 0 \\ 0 & 1
\end{pmatrix}, \pm\begin{pmatrix}
    1 & 1 \\ 1 & 2
\end{pmatrix},\pm \begin{pmatrix}
    0 & 2 \\ 1 & 0
\end{pmatrix},\pm\begin{pmatrix}
    1 & 2 \\ 2 & 2
\end{pmatrix}\right\}.
\]
From \cite[exercise 12.15 (d)]{Cox1997}, we have
$$
\Gamma_{\gamma_2} = \pi_3^{-1}(G_{\gamma_2}) \cap \SL_2(\Z).
$$

\begin{lemma}\label{lem: gamma 2 normalizer}
    The matrix $D(1,2)$ lies in the normalizer of $G_{\gamma_2}$. The subgroup $H = \langle D(1,2), G_{\gamma_2} \rangle \subseteq \GL_2(\Z/3\Z)$ satisfies $H\cap \SL_2(\Z/3\Z) = G_{\gamma_2}$ and $\det H = (\Z/3\Z)^\times$.
\end{lemma}
\begin{proof}
    Given a matrix $g=\begin{pmatrix}
        a & b \\ c & d
    \end{pmatrix}$, we have
    $$
    D(1,2)gD(1,2)^{-1}=\begin{pmatrix}
        a & 2b \\ 2c & d
    \end{pmatrix}.
    $$
    Taking~$g$ to be one of the above generators of $G_{\gamma_2}$, we see that $D(1,2)$ indeed lies in the normalizer of $G_{\gamma_2}$. Thus any element~$h\in H$ can be written in the form $D(1,2)^r g$ where~$r\in \{0,1\}$ and~$g\in G_{\gamma_2}$. Since $\det D(1,2) = 2$, this shows $H\cap \SL_2(\Z/3\Z) = G_{\gamma_2}$ and $\det H = (\Z/3\Z)^\times$.
\end{proof}

\begin{lemma}\label{lem: gamma 2 morphism}
Let~$H$ be as in \cref{lem: gamma 2 normalizer}.
    The cover $X_H^{(-1)}\to X(1)$ has degree 3, and $\gamma_2$ defines an isomorphism $X_{H}^{(-1)}\to \Proj^1$ over $\Q$.
\end{lemma}
\begin{proof}
    The degree of the cover $X_{\Gamma_{\gamma_2}}^{(-1)}\to X(1)$ is equal to
    $$[\PSL_2(\Z):(\Gamma_{\gamma_2} \mod \pm I)] = [\SL_2(\Z): \Gamma_{\gamma_2}] = [\SL_2(\Z/3\Z):G_{\gamma
    _2}] = \frac{24}{8} = 3.
    $$
    Since we have the relation $j = \gamma_2^3$, \cref{prop: functions GAGA} \ref{item: Hauptmodul} shows that $\gamma_2$ is an isomorphism.
    Recall that
    \begin{equation}\label{eq: q exp E4}
    E_4(\tau)=\frac{3g_2(\tau)}{4\pi^4}
    =1+240\sum_{n\geq1}\sigma_3(n)q^n
    \in\Z[[q]].
    \end{equation}
    Hence, by the product expansion of $\eta$,
    \begin{equation}\label{eq: q exp gamma 2}
    \gamma_2(\tau)
    =\frac{E_4(\tau)}{\eta(\tau)^8}
    =q^{-1/3}E_4(\tau)
    \prod_{n\geq1}(1-q^n)^{-8}
    \in\mathbb Z((q^{1/3})).
    \end{equation}
    This means that $\gamma_2$ is defined over $\Q$, by \cref{prop: functions GAGA} \ref{item: q expansion principle}.
\end{proof}
\begin{proposition}
Taking $H = \langle D(1,2), G_{\gamma_2} \rangle$, we can define the polynomial $\Phi_\ell^{\gamma_2}$ for any prime~$\ell\neq 3$. Moreover, this polynomial has coefficients in $\Z$, and it is denominator-free and symmetric.
\end{proposition}
\begin{proof} We observe that 
$$D(2,1)=\pi_3(S)D(1,2)\pi_3(S)^{-1}\in H, \quad \text{so}\quad 2I=D(2,1)D(1,2)\in H.$$
Taking $M_u=D(1,u)$ where $u = (\ell\mod 3)$ and $v=w=0$, we find that $H$ satisfies both conditions of \cref{thm: summary conditions}. Therefore, $\Phi_\ell^{\gamma_2}$ is denominator-free and symmetric.  In addition, by \eqref{eq: q exp E4} and \eqref{eq: q exp gamma 2}, the $q$-expansion of $\gamma_2$ has integer coefficients and leading coefficient 1; hence, by \cref{cor: coefficients in Z}, $\Phi_\ell^{\gamma_2}$ has coefficients in $\Z$.
\end{proof}

When $\ell=2$, Elkies computed in \cite[equation (22)]{Elkies1997EllipticAM}:
$$
\Phi_2^{\gamma_2}(x,y)= x^3+y^3-(xy)^2+495xy-\numprint{54000}.
$$
Since the cover $X_{H}^{(-1)}\to X(1)$ has degree 3, applying \cref{thm:main-height} to $\gamma_2$, we then obtain a bound that is asymptotically slightly looser than \cite[Theorem 1.4]{zhang2026explicitheightboundsmodular}, as Zhang's second term is $O(\ell)$ while ours is $O(\ell \sqrt{\log \ell})$.

\subsection{The Weber modular function}

Define the functions
$$ \f(\tau) := \zeta_{48}^{-1}\frac{\eta((\tau+1)/2)}{\eta(\tau)},\qquad  \f_1(\tau) := \frac{\eta(\tau/2)}{\eta(\tau)},\qquad\f_2(\tau) := \sqrt{2}\frac{\eta(2\tau)}{\eta(\tau)}.
$$
The function $\f$ is called the \textit{Weber modular function}. In this section, we apply our results from §\ref{section: construction} to the function $\f$ and the associated modular polynomial $\Phi_\ell^\f$.

Define an action of $\SL_2(\Z)$ on the vector space of meromorphic functions by $(\gamma, f)\mapsto f\circ \gamma^{-1}$.
Cox proves in \cite[Cor. 12.19]{Cox1997} that $\SL_2(\Z)$ acts on the vector space $\C \f \oplus \C \f_1 \oplus \C \f_2 $ as 
\begin{equation}\label{representation Weber}
\rho: T\mapsto \begin{pmatrix}
    0 & \zeta_{48} & 0 \\ \zeta_{48} & 0 & 0 \\ 0 & 0 & \zeta_{48}^{-2}
\end{pmatrix},~~~S\mapsto  
\begin{pmatrix}
    1 & 0 & 0 \\ 0 & 0 & 1 \\ 0 & 1 & 0
\end{pmatrix},\end{equation}
where $\zeta_{48} = e^{2\pi i /48}$ and $\rho:\SL_2(\Z)\to \GL_3(\Q(\zeta_{48}))$.
The largest congruence subgroup $\Gamma_\f$ leaving $\f$ invariant is then the group of matrices fixing $(1, 0, 0)^t$. In other words,
$$\Gamma_\f=\rho^{-1}(\Sg_1),\qquad \text{ where }\qquad\Sg_1:= \left\{\begin{pmatrix}
    1 & * & * \\ 0 & * & * \\ 0 & * & *
\end{pmatrix}\in \GL_3(\Q(\zeta_{48}))\right\}.$$
\begin{lemma}\label{lem: ker representation}
    We have $\Gamma(48)\subseteq \ker \rho$. In particular, $\Gamma(48)\subseteq \Gamma_\f$ and $\f$ is a modular function of level $48$.
\end{lemma}
\begin{proof}
    Using the notation from \cite{Enge.Morain2014}, we have $\f_1=\m_2$. By \cite[theorem 3.2]{Enge.Morain2014}, taking $N=2$, $e=1$, $t=24$, $s=24$, $\m_2$ is modular for $\Gamma(24)\cap \Gamma^0(48)$, and in particular for $\Gamma(48)$. Using \cref{representation Weber}, we have $\f = T^{-1}\zeta_{48}\f_1$ and $\f_2=S^{-1}\f_1$. We then have
    $$\Gamma_{\f}=T^{-1}\Gamma_{\f_1}T \qquad\text{and}\qquad\Gamma_{\f_2}=S^{-1}\Gamma_{\f_1}S.
    $$
    Since $\Gamma(48)$ is normal in $\SL_2(\Z)$, this implies that $\Gamma(48)\subseteq \Gamma_{\f}$ and $\Gamma(48)\subseteq \Gamma_{\f_2}$.
\end{proof}
Given \cref{lem: ker representation}, the representation $\rho:\SL_2(\Z) \to \GL_3(\Q(\zeta_{48}))$ factors through $\SL_2(\Z/48\Z)\simeq \SL_2(\Z)/\Gamma(48)$. 
\begin{proposition}\label{prop: largest degree}
    The degree of the cover $X_{\Gamma_\f}\to X(1)$ is equal to $72$ and is maximal among genus 0 modular covers of $X(1)$.
\end{proposition}
\begin{proof}
       Let $G_\f:= \pi_{48}(\Gamma_\f)$. Since $G_\f=\Stab_{\SL_2(\Z/48\Z)}(\f)$, we have
    $$
    [\PSL_2(\Z):(\Gamma_\f \mod \pm I)]=[\SL_2(\Z):\Gamma_\f]=[\SL_2(\Z/48\Z):G_\f]=|\SL_2(\Z/48\Z)\cdot\f|.
    $$
    The first equality holds since $-I\in \Gamma_\f$, the second by reduction modulo 48, and the third by the orbit-stabilizer theorem.
    To compute the last term, note that for any $k\in \Z$, we have $T^{2k}\f=\zeta_{48}^{2k}\f$. Therefore, the orbit of $\f$ contains the elements 
    \begin{equation}
    \zeta_{48}^{2k}\f, \qquad \text{for }k\in \Z/24\Z.
    \label{eq: orbit f 1}
    \end{equation}
    Applying $T$ and subsequently $S$ to these elements gives the terms
    \begin{equation}\label{eq: orbit f 2}
    \zeta_{48}^{2k+1}\f_1, \qquad\zeta_{48}^{2k+1} \f_2
    \end{equation}
    in the orbit of $\f$. Yet, we note that the set formed by the elements in \cref{eq: orbit f 1} and~\cref{eq: orbit f 2} is stable under $S$ and $T$. Therefore, it constitutes the orbit of $\f$, and contains $3\cdot24 = 72$ elements.
    By the classification from \cite{CHUA2004408}, any genus 0 modular cover of $X(1)$ must have degree smaller than or equal to 72.
\end{proof}

Our next step is to construct a suitable subgroup $H\leq \GL_2(\Z/48\Z)$ such that $H \cap \SL_2(\Z/48\Z) = \pi_{48}(\Gamma_\f)$.

\begin{proposition}\label{Cor: Weber normalizer}
    For any $u\in (\Z/48\Z)^\times$, we have $$D(1,u)\in N_{\GL_2(\Z/48\Z)}(\pi_{48}(\Gamma_\f)).$$
    In particular, the set
    \begin{equation}\label{eq: H Weber}
    H= \{D(1,u):u\in (\Z/48\Z)^\times\}\pi_{48}(\Gamma_\f)
    \end{equation}
    is a group satisfying
    \begin{equation}\label{eq: H Weber ppties}
    \Gamma_H=\Gamma_\f \quad\text{and}\quad(\Z/48\Z)^\times\subseteq \det H.
    \end{equation}
\end{proposition}

The proof will rely on two lemmas.

\begin{lemma}
\label{lem:epsilon-u}
    For every $u\in (\Z/48\Z)^\times$, we have $u^2\in \{1,25\}$. In particular, $\zeta_{48}^{u^{-1}-u}=\epsilon_u$, where
$$\epsilon_u=\begin{cases}1 &\text{ if }u^2\equiv 1 \mod 48 \\ -1 &\text{ if } u^2 \equiv 25 \mod 48.\end{cases}$$

\end{lemma}
\begin{proof}
    Since $3\nmid u$, we have $u^2\equiv1\mod3$. On the other hand, every odd square is congruent to $1$ or $9$ modulo $16$; the assertion follows from the Chinese remainder theorem.
    Consequently, for any $u\in (\Z/48\Z)^\times$, we have $$u^{-1}-u=u^{-1}(1-u^2)= \begin{cases}
        0 &\text{ if }u^2\equiv 1 \mod 48 \\ 24u^{-1} &\text{ if }u^2\equiv 25 \mod 48.
    \end{cases}$$
    Since $u^{-1}$ is odd, $(\zeta_{48}^{24})^{u^{-1}}=-1$, so $\zeta_{48}^{u^{-1}-u}=\epsilon_u.$
\end{proof}

For $u\in (\Z/48\Z)^\times$,
let $\sigma_u\in \Gal(\Q(\zeta_{48})/\Q)$ be defined by $\sigma_u: \zeta_{48} \mapsto \zeta_{48}^u$, and define
$$C_u = \begin{pmatrix}
    1 & 0 & 0 \\ 0 & 1 & 0 \\ 0 & 0 & \epsilon_u
\end{pmatrix}\in \GL_3(\Q(\zeta_{48}))$$
where~$\epsilon_u$ is as in \cref{lem:epsilon-u}.

\begin{lemma}
\label{lem:weber-conj}
    For any $u\in (\Z/48\Z)^\times$ and
    $g\in \SL_2(\Z/48\Z)$, define
    $$c_u(g):=D(1,u)gD(1,u)^{-1}.$$
    Then
    $$
    \rho(c_u(g))
    =C_u\sigma_{u^{-1}}(\rho(g))C_u^{-1}.
    $$
\end{lemma}

\begin{proof}
    Let $u'=u^{-1}$.
    It suffices to prove the identity for $S$ and $T$.
    We have
    $$
    c_u(S)=T^{u'}ST^uST^{u'}.
    $$
    Since $u$ and $u'$ are odd,
    $$
    \rho(T^u)=
    \begin{pmatrix}
        0 & \zeta_{48}^u & 0\\
        \zeta_{48}^u & 0 & 0\\
        0 & 0 & \zeta_{48}^{-2u}
    \end{pmatrix},
    \qquad
    \rho(T^{u'})=
    \begin{pmatrix}
        0 & \zeta_{48}^{u'} & 0\\
        \zeta_{48}^{u'} & 0 & 0\\
        0 & 0 & \zeta_{48}^{-2u'}
    \end{pmatrix}.
    $$
    Therefore,
    $$
    \rho(c_u(S))
    =
    \begin{pmatrix}
        \zeta_{48}^{2(u'-u)} & 0 & 0\\
        0 & 0 & \zeta_{48}^{u-u'}\\
        0 & \zeta_{48}^{u-u'} & 0
    \end{pmatrix}
    =
    \begin{pmatrix}
        1 & 0 & 0\\
        0 & 0 & \epsilon_u\\
        0 & \epsilon_u & 0
    \end{pmatrix}=C_u\sigma_{u'}(\rho(S))C_u^{-1},
    $$
    where we used $u\equiv u'\pmod{24}$ and
    $\epsilon_u=\zeta_{48}^{u'-u}\in\{\pm1\}$.
    On the other hand, $$c_u(T)=T^{u'},$$
    so
    $$
    \rho(T^{u'})
    =\sigma_{u'}(\rho(T))
    =C_u\sigma_{u'}(\rho(T))C_u^{-1},
    $$
    since $C_u$ commutes with $\rho(T^{u'})$.
\end{proof}

\begin{proof}[Proof of \cref{Cor: Weber normalizer}]
If $g\in \pi_{48}(\Gamma_\f)$, then $\rho(g) \in \Sg_1$. Therefore, for any $u\in (\Z/48\Z)^\times$, we have $\sigma_{u^{-1}}(\rho(g))\in \Sg_1$, and in particular $C_u \sigma_{u^{-1}}(\rho(g)) C_u^{-1}\in \Sg_1$, meaning that $\rho(c_u(g))\in \Sg_1$ and $c_u(g)\in \pi_{48}(\Gamma_\f)$. This proves the first assertion of \cref{Cor: Weber normalizer}; the two others are easy consequences.
\end{proof}

Once~$H$ has been constructed, the properties of the Weber modular polynomials are easily deduced from our general statements.

\begin{proposition}
Let~$H$ be as in \cref{Cor: Weber normalizer}.
    The Weber modular function~$\f$ defines an algebraic isomorphism $X_{H}^{(-1)} \to \Proj^1$ defined over $\Q$.
\end{proposition}
\begin{proof}
    From \cite[Theorem 12.7]{Cox1997}, we have 
    $$j=\frac{(\f^{24}-16)^3}{\f^{24}}.$$
    In particular, the relation between $\f$ and $j$ has degree $72=\deg(X_H^{(-1)}\to X(1))$.
    In addition, the product expansion of $\eta$ yields 
    \begin{equation}\label{eq: q exp Weber}
    \f(\tau)
    =q^{-1/48}\prod_{n\geq1}(1+q^{n-1/2})
    \in\Z((q^{1/48})).
    \end{equation}
    Therefore, by \cref{prop: functions GAGA} \ref{item: q expansion principle} and \ref{item: Hauptmodul}, $\f$ is an isomorphism $X_{H}^{(-1)}\to \Proj^1_\Q$.   
\end{proof}

\begin{proposition}
    Taking~$H$ as in \cref{Cor: Weber normalizer}, we can define the polynomial $\Phi_\ell^{\f}$ for any prime $\ell\notin\{2,3\}$. This modular polynomial has coefficients in $\Z$ and it is denominator-free and symmetric.
\end{proposition}
\begin{proof}
    Since $\pi_{48}(S)\in \pi_{48}(\Gamma_\f)$, we also have $$D(u,1)=\pi_{48}(S)D(1,u)\pi_{48}(S)^{-1}\in H.$$
    In particular, $$D(u^{-1},u)=D(u,1)^{-1}D(1,u)\in H.$$ 
    Taking $v=w=0$ and using $D(1,u)$ and $D(u^{-1},u)$, by \cref{equivalence}, $\Phi_\ell^\f$ is denominator-free. On the other hand, $uI= D(u,1)D(1,u)\in H$, which proves that $\Phi_\ell^\f$ is symmetric. In addition, by \cref{eq: q exp Weber}, the $q$-expansion of $\f$ has integer coefficients and leading coefficient 1; hence, by \cref{cor: coefficients in Z}, $\Phi_\ell^{\f}$ has coefficients in $\Z$.
\end{proof}

The values of $\Phi_\ell^\f$ for $5\leq \ell\leq 1000$ have been computed in \cite{sutherland_weber_modular_polynomials}. For example,
$$
\Phi_5^f(x,y) = x^6-x^5y^5+xy+y^6.
$$
\noindent Applying \cref{thm:main-height} to $\Phi_\ell^\f$, we again find a slightly looser bound than the one in \cite[Theorem 1.6]{zhang2026explicitheightboundsmodular}, since our second term involves $O(\ell \sqrt{\log \ell})$ while Zhang's involves $O(\ell)$.

\printbibliography[heading=bibintoc,
title={References}]

\appendix
\section{Modular curves as moduli spaces}
\label{appendix}

\markboth{APPENDIX}{APPENDIX}
\markright{APPENDIX}

\subsection{Level structures}
\label{subsec:appendix-level-structures}

In order to define modular curves as moduli spaces, we need to define elliptic curves and their level $H$ structures, where $H$ is a subgroup of $\GL_2(\Z/n\Z)$.
This theory is developed in the foundational works \cite{Deligne.Rapoport1973,KatzMazur} over $\Z[1/n]$-schemes. In this appendix, we work with $\Q$-schemes throughout for simplicity.

An \textit{elliptic curve} over a $\Q$-scheme $S$ is a proper smooth morphism of schemes $E\to S$ with geometrically connected fibers of genus 1, together with a section $0:S\to E$.  This defines a unique structure of abelian group scheme for
which 0 is the neutral element \cite[(2.1.2)]{KatzMazur}. Therefore, we have a multiplication-by-$n$ $S$-homomorphism $[n]:E\to E$. It is finite and locally free of rank $n^2$. The kernel $E[n]$ of $[n]$ is finite étale over $S$, and is locally isomorphic to $(\Z/n\Z)_S^2:=\coprod_{(a,b)\in (\Z/n\Z)^2}S$ for the étale topology on $S$  \cite[(2.3.1)]{KatzMazur}.

A \textit{full level $n$ structure}, if it exists, is then a global isomorphism $$E[n]\to (\Z/n\Z)_S^2$$ 
of finite étale $S$-group schemes \cite[(IV. 2.3)]{Deligne.Rapoport1973}.

Level structure for subgroups $H\leq \GL_2(\Z/n\Z)$ are defined as follows. Let $\Sch_S$ be the category of schemes over $S$, and consider the sheaf of sets $\F_n(E)$ on $\Sch_S$ defined by
$$
T\mapsto \{\text{full level $n$ structures on }E_T\},
$$
where $E_T=E\times_ST$.
Any matrix $h\in\GL_2(\Z/n\Z)$ induces an isomorphism $(\Z/n\Z)_T^2\to (\Z/n\Z)_T^2$, also denoted by~$h$, by permuting the components $T$ indexed by elements of $(\Z/n\Z)^2$ considered as column vectors. This induces a right action of $\GL_2(\Z/n\Z)$ on $\F_n(E)$ given by
\begin{equation}
    \begin{matrix}
        \F_n(E)(T) \times \GL_2(\Z/n\Z) & \to & \F_n(E)(T) \\
        (\alpha,h) &\mapsto &h^t\circ \alpha.
    \end{matrix}
\label{eq:H-action}
\end{equation}
We then consider the étale sheafification $\F_H(E)$ of the presheaf of sets 
$$T\mapsto \F_n(E)(T)/H.$$
A \textit{level~$H$ structure} is defined as an element of $\F_H(E)(S)$; this differs from \cite[IV. 3.1]{Deligne.Rapoport1973} by a transpose. The special case $H = \{1\}$ corresponds to full level~$n$ structures as above.

Two elliptic curves $E,E'$ over~$S$ endowed with level~$H$ structures $\alpha,\alpha'$ are declared \emph{isomorphic} if there exists an isomorphism $\phi:E\to E'$ such that, under the induced isomorphism $\phi^*: \F_H(E')\to \F_H(E)$, we have $\phi^*\alpha' = \alpha$.

By the properties of sheafification, level $H$ structures can be described as follows.
Given schemes $S''\to S'\to S$ and a section $\alpha \in \F_n(E)(S')$, denote by $\alpha|_{S''}=\alpha \times_{S'} S''$.
Consider pairs $((S_i\to S)_i, (\alpha_i)_i)$, where $(S_i\to S)_i$ is an étale cover and $(\alpha_i)_i\in \prod_i \F_n(E)(S_i)$. Two such pairs $((S_i\to S)_i, (\alpha_i)_i)$, $((S_j'\to S)_j, (\alpha_j')_j)$ are said to be $H$-equivalent if for every $i,j$, there exists a fine enough étale cover $(U_k\to S_i\times_S S_j')_k$ with elements $(h_k)_k\subseteq H$ such that
$$
(\alpha_i|_{S_i\times_S S_j'})|_{U_k}=h_k^t \circ (\alpha_j'|_{S_i\times_S S_j'})|_{U_k}.
$$
A pair $((S_i\to S)_i, (\alpha_i)_i)$ is \textit{admissible} if it is $H$-equivalent to itself. Then a level~$H$ structure exactly corresponds to the data of an admissible pair $((S_i\to S)_i, (\alpha_i)_i)$ up to $H$-equivalence.

This general definition indeed corresponds to the definition of level~$H$ structures given in §\ref{subsec: level struct for fields} over fields:

\begin{proposition}
\label{prop:galois-sheaf}
        For every field $k$ containing $\Q$, we have
        $$
        \F_H(E)(k) \simeq (\F_n(E)(\overline k)/H)^{G_k},
        $$
        where $G_k = \Gal(\overline k/k)$. In other words, $\F_H(E)(k)$ is the set of $H$-classes of isomorphisms $\alpha:E[n]_{\overline k}\to (\Z/n\Z)^2_{\overline k}$ such that for any $\sigma \in G_k$ there exists $h_\sigma \in H$ such that
        $$
        \alpha\circ \rho_{E,n}(\sigma) = h_\sigma^t \circ \alpha.
        $$
    \end{proposition}
    \begin{proof}
        We can work over a finite Galois extension $L/k$ satisfying $E[n](L)=E[n](\overline k)$, which exists since $E[n]$ is finite étale. In particular, we have
        $$
        E[n]_L \simeq (\Z/n\Z)^2_L,
        $$
        meaning that $\F_n(E)$ is constant on $\Sch_L$. Consequently, we have
        $$
        \F_H(E)(L) \simeq \F_n(E)(L)/H.
        $$
        Consider the finite étale cover $\Spec L \to \Spec k$. Then, as $\F_H(E)$ is a sheaf, we have the equalizer
        \begin{equation}
        \F_H(E)(k) \to \F_H(E)(L) \rightrightarrows \F_H(E)(L \otimes_k L).
        \label{eq: equalizer diragram}
        \end{equation}
        Moreover, since $L/k$ is Galois, we have 
        $$
        L\otimes_k L \simeq \prod_{\sigma \in \Gal(L/k)} L.
        $$
        Under this identification, diagram \cref{eq: equalizer diragram} describes the elements of $\F_H(E)(k)$ as the elements $[\alpha]\in \F_H(E)(L)$ satisfying $[\alpha \circ \rho_{E,n}(\sigma)]=[\alpha]$ for all $\sigma\in \Gal(L/k)$, which gives the result.
    \end{proof}

\subsection{The Weil pairing}
\label{subsec:appendix-weil}

As in §\ref{subsec: Weil pairing}, we use the Weil pairing to distinguish between connected components of modular curves. For any~$\Q$-scheme $S$ and elliptic curve~$E/S$, the Weil pairing is a morphism $e_n: E[n]\times E[n] \to \mu_{n,S}$, where $\mu_{n,S}$ denotes the kernel of the multiplication by $n$ in $\mathbb{G}_m$, the multiplicative group scheme  over $S$; it is constructed in \cite[(2.8.5.1)]{KatzMazur}. In the rest of this subsection, we assume that the base scheme~$S$ is a $K_H$-scheme, where~$K_H$ is constructed as in~§\ref{subsec: Weil pairing}.

Let~$\mu_{n,\,\Spec K_H}^0 \subseteq \mu_{n,\,\Spec K_H}$ denote the subscheme of primitive $n$-th roots of unity. As a set of $\overline \Q$-points equipped with a Galois action, this scheme~$\mu_{n,\,\Spec K_H}^0$ is identified with $\{\zeta_n^\beta: \beta\in (\Z/n\Z)^\times\}$ equipped with the multiplicative action of $\Gal(\Q(\zeta_n)/K_H))\simeq \det H$. We denote the quotient scheme by
$$
    \pi_{0,H} = \mu_{n,\,\Spec K_H}^0 / \det H \simeq \bigsqcup_{\beta\in (\Z/n\Z)^\times/\det H} \Spec K_H.
$$

If~$S$ is a $K_H$-scheme and $\alpha$ is a level~$H$ structure on some elliptic curve $E/S$ encoded by $((S_i\to S)_i, (\alpha_i)_i)$, we construct a morphism $e_n(\alpha): S\to \pi_{0,H}$ as follows. For each~$i$, we have an isomorphism $\alpha_i: E_{S_i}[n] \to (\Z/n\Z)_{S_i}$, so we have two sections $P_1,P_2: S_i\to E_{S_i}[n]$ given by $\alpha_i^{-1}(1,0)^t$ and $\alpha_i^{-1}(0,1)^t$ respectively. Then $e_n(P_1,P_2)$ yields a morphism $e_n(\alpha)_i: S_i\to \mu^0_{n,S_i}$, as the Weil pairing is nondegenerate \cite[(3.20.1)]{Deligne.Rapoport1973}.

We compose this morphism with the structural morphism $\mu^0_{n,S_i} \to \mu^0_{n,\,\Spec K_H}$ and the quotient map $\mu^0_{n,\,\Spec K_H}\to \pi_{0,H}$. The admissibility condition guarantees that the resulting morphisms $e_n(\alpha)_i$ glue into a morphism $e_n(\alpha): S\to \pi_{0,H}$, and $K_H$ was precisely constructed so that $e_n(\alpha)$ does not depend on the choice of $H$-equivalent representative of~$\alpha$ \cite[(3.20.3)]{Deligne.Rapoport1973}.  Further, if~$\phi:(E,\alpha)\to (E',\alpha')$ is an isomorphism between elliptic curves with level~$H$ structures over~$S$, then $e_n(\alpha) = e_n(\alpha')$.

Given the above description of~$\pi_{0,H}$, for any~$\beta\in (\Z/n\Z)^\times/\det H$, we abuse notation and write ``$e_n(\alpha) = \zeta_n^\beta$'' if~$e_n(\alpha)$ takes values in the connected component of $\pi_{0,H}$ corresponding to $\beta$. This is consistent with our usage in §\ref{subsec:modular-fields}.
As above, this definition can be simplified if~$S = \Spec k$, where~$k$ is a field containing $K_H$, and it agrees with our construction in §\ref{subsec: Weil pairing}.

\subsection{Fine and coarse moduli spaces}
\label{subsec:appendix-moduli}

For the reader's convenience, we recall general concepts pertaining to moduli spaces before applying them to the case of modular curves.

Let $S_0$ be a scheme. A \textit{moduli problem} over $S_0$ is a contravariant functor $\Pa:\Sch_{S_0}\to \Sets$, where $\Sch_{S_0}$ is the category of schemes over $S_0$. A moduli problem should be seen as a problem consisting of classifying sets of objects over varying $S_0$-schemes. 

A \textit{fine moduli space} for $\Pa$ is a scheme $M_0\in \Sch_{S_0}$ such that $\Pa \simeq \Hom_{S_0}(-,M_0)$. In other words, for any $S\in \Sch_{S_0}$, we have a bijection $\Pa(S)\to M_0(S)$ that is moreover natural in $S$.
In particular, $\Pa(M_0)$ contains a universal element corresponding to the identity map on~$M_0$. Elements of $\Pa(S)$ are then obtained uniquely from base change of the universal element along morphisms $S\to M_0$. Intuitively, $M_0$
``solves'' $\Pa$ in a strong sense.

It is well-known that fine moduli spaces do not exist for every moduli problem, see for example \cite[§1.8]{Kollár_2023}, \cite[§23]{Hartshorne2010}, but a weaker notion exists. A \textit{coarse moduli space} for~$\Pa$ is an $S_0$-scheme $M$ endowed with a natural transformation $t:\Pa \to \Hom_{S_0}(-,M)$ such that:
\begin{enumerate}
    \item Whenever $S=\Spec k$ for an algebraically closed field $k$, the map $t_S:\Pa(S) \to M(S)$ is a bijection.
    \item \label{factor coarse moduli space} For any other $S_0$-scheme $M'$ and natural transformation $t':\Pa \to M'$, there exists a unique natural transformation $\eta:\Hom_{S_0}(-,M) \to \Hom_{S_0}(-,M')$ such that $\eta \circ t = t'$.
\end{enumerate}
A fine moduli space, if it exists, is also a coarse one; and a coarse moduli space, if it exists, is unique up to unique isomorphism by the Yoneda lemma.

\begin{remark}
\label{rem:coarse-morphisms}
    If we have two moduli problems $\Pa,\Pa'$ with coarse moduli spaces given by $M,t$ and $M',t'$ respectively, and a natural transformation $\mathcal{T}:\Pa \to \Pa'$, then by the property \ref{factor coarse moduli space}, the composition $\Pa \to \Pa' \to \Hom_{S_0}(-,M')$ factors through $\Hom_{S_0}(-,M)$, yielding a natural transformation $\Hom_{S_0}(-,M)\to \Hom_{S_0}(-,M')$. By the Yoneda lemma, this gives an $S_0$-scheme morphism $M\to M'$. On $S$-points, the induced map $M(S)\to M'(S)$ then corresponds to the map $\Pa(S)\to \Pa'(S)$.
\end{remark}

\subsection{Modular curves}
\label{subsec:modular-curves}

Modular curves arise as coarse moduli spaces of elliptic curves endowed with level structures. Following \cite[(7.1.3)]{KatzMazur}, given a subgroup $H \leq \GL_2(\Z/n\Z)$, we consider the moduli problem over $\Sch_{\Q}$ given by
\begin{align*}
\Pa_H: S \mapsto \{&\text{isomorphism classes of pairs }(E,\alpha): \\
&\text{$E$ is an elliptic curve over $S$, $\alpha \in \F_H(E)(S)$}\}.
\end{align*}
We are also interested in classifying level structures with fixed Weil pairing. Following \cite[(9.4.2)]{KatzMazur}, we work over $S_0 = \Spec K_H$. For $\beta \in (\Z/n\Z)^\times/\det H$, we define the following moduli problem over~$\Sch_{S_0}$:
\begin{align*}
\Pa_H^{(\beta)}: S \mapsto \{&(E,\alpha)\in \Pa_H(S): e_n(\alpha) = \zeta_n^\beta\}.
\end{align*}

The properties of modular curves used in this paper can then be summarized as follows.

\begin{theorem}
    \begin{enumerate}
        \item The moduli problem $\Pa_H$ admits a coarse moduli space $Y_H$, which is a smooth algebraic curve over~$\Q$. For any $\beta\in (\Z/n\Z)^\times/\det H$, the moduli problem $\Pa_H^{(\beta)}$ admits a coarse moduli space $Y_H^{(\beta)}$, which is a smooth algebraic curve over $K_H$.
        \item Let~$X_H$ and $X_H^{(\beta)}$ denote the smooth projective closures of $Y_H$ and~$Y_H^{(\beta)}$ respectively. Then~$X_H$ is connected; its geometric connected components are individually defined over the number field $K_H$, and are exactly the $X_H^{(\beta)}$ for $\beta\in (\Z/n\Z)^\times/\det H$.
    \end{enumerate}
\end{theorem}

\begin{proof}
    See \cite[(7.1.2), (8.1.1), (8.1.5)]{KatzMazur} for properties of~$X_H$ and \cite[(9.1.10), (9.4.2), (9.4.3.1)]{KatzMazur} for properties of $X_H^{(\beta)}$.
\end{proof}

\subsection{Modular curves and isogenies}

Finally, we explain why the diagram \cref{eq:diagram-YH} indeed arises from a morphism of algebraic curves over $K_H$. To this end, we use the moduli interpretation and \cref{rem:coarse-morphisms}: we will construct two natural transformations $\mathcal{P}_{H'} \to \mathcal{P}_H$ corresponding to the domain and codomain of $\ell$-isogenies between elliptic curves (over $\Q$-schemes) equipped with level~$H$ structures.

We first review isogenies in this setting. Let~$E,E'$ be elliptic curves over a base scheme~$S$. An \emph{$\ell$-isogeny} $\phi:E\to E'$ is a finite, locally free, surjective homomorphism $E\to E'$ over~$S$ such that $\phi_* \OO_E$ has rank~$\ell$ as an $\OO_{E'}$-module \cite[(2.8.1) p.\,87]{KatzMazur}. This notion agrees with the usual one over fields. There is an equivalence between isogenies and subgroup schemes:
\begin{itemize}
    \item If~$\phi: E\to E'$ is an $\ell$-isogeny, then $\ker\phi\subseteq E[\ell]$ is a finite locally free $S$-subgroup scheme of~$E$ of order~$\ell$ \cite[proof of Thm.~2.5.1 p.\,80 and Lemma 1.4.2]{KatzMazur}.
    \item Conversely, if~$C\subseteq E$ is a finite locally free $S$-subgroup scheme of order~$\ell$, then the quotient $E' = E/C$ is an elliptic curve over~$S$, and the quotient map $E\to E'$ is an $\ell$-isogeny; moreover the formation of~$E/C$ commutes with arbitrary base change. See \cite[SGA3, Exp.~V, Thm.~4.1]{artinSeminairesGeometrieAlgebrique1971} and \cite[(3.3) p.\,100]{KatzMazur}.
\end{itemize}

Under our running assumption $\ell\nmid n$, any $\ell$-isogeny $\phi:E\to E'$ induces an isomorphism~$\phi^*$ between $E[n]$ and~$E'[n]$, hence an isomorphism between $\F_H(E')$ and $\F_H(E)$. If~$\alpha$ and $\alpha'$ are level~$H$ structures on~$E$ and~$E'$ respectively, we say that $\phi$ \emph{respects}, or \emph{is compatible with}, these structures if~$\phi^*\alpha' = \alpha$. If~$\alpha$ is a fixed level~$H$ structure on~$E$, then there exists a unique level~$H$ structure~$\alpha'$ on~$E'$ such that $\phi$ is compatible with~$\alpha$ and~$\alpha'$. Combining this remark with the above equivalence between isogenies and subgroup schemes justifies the following definition.

\begin{definition}
    Let~$H$ be a subgroup of~$\GL_2(\Z/n\Z)$, let~$\ell\nmid n$ be a prime, and let $S_0$ be a base scheme. We define the moduli problem $\Pa_H^{\ell\mathrm{-isog}}$ on $\mathrm{Sch}_{S_0}$ as follows:
\begin{align*}
\Pa_H^{\ell\mathrm{-isog}}: S \mapsto \{&\text{isomorphism classes of triples }(E,\alpha,C): \\
&\text{$E$ is an elliptic curve over $S$, $\alpha \in \F_H(E)(S)$}, \\
&\text{and $C$ is a finite locally free $S$-subgroup scheme of~$E$ of order~$\ell$}\}.
\end{align*}
\end{definition}

We then have two natural transformations $\Pa_H^{\ell\mathrm{-isog}} \to \Pa_H$ given by the domain and codomain maps. We may describe these natural transformations by maps (between sets) $\Pa_H^{\ell\mathrm{-isog}}(S) \to \Pa_H(S)$ for any $S_0$-scheme~$S$, as follows:
\begin{equation}
    \begin{tikzcd}
        & \Pa_H^{\ell\mathrm{-isog}}(S) \ar[ld, swap, "{(E,\alpha,C) \mapsto (E,\alpha)}"] 
        \ar[rd, "{(E,\alpha,C)\mapsto (E/C, \alpha')}"] & \\
        \Pa_H(S) && \Pa_H(S)
    \end{tikzcd}
\label{eq:moduli-diagram}
\end{equation}
where, on the right, $\alpha'$ is the unique level~$H$ structure on~$E/C$ such that the quotient map $E\to E/C$ is compatible with $\alpha,\alpha'$. The maps in~\eqref{eq:moduli-diagram} are indeed compatible with base changes $S'\to S$ thanks to the above compatibility between base change and quotients, and correspond to the maps from \cref{eq:diagram-YH} and \cref{prop: maps over C} over~$\C$.

In order to translate the diagram~\eqref{eq:moduli-diagram} into morphisms between the modular curves $Y_{H'}$ and $Y_H$, where~$H'$ is taken as in \cref{eq: def H'}, we need to find an isomorphism between $\Pa_H^{\ell\mathrm{-isog}}$ and~$\Pa_{H'}$.

\begin{proposition}
\label{prop:moduli-product}
    Let~$S$ be a $\Q$-scheme, and let~$(E,\alpha)$ be an elliptic curve over~$S$ with level~$H'$ structure, where~$\alpha$ is represented by an admissible pair $((S_i\to S)_i, (\alpha_i)_i)$. For each~$i$, consider the isomorphism
    \[
    \begin{matrix}
    (\iota_n, \iota_\ell):& E_{S_i}[n\ell] & \to & E_{S_i}[n] \times E_{S_i}[\ell]\\
   & P & \mapsto & (\ell P, nP).
    \end{matrix}
    \]
    For~$i\in I$, let~$\alpha'_i = \pi_n \circ \alpha_i \circ \iota_n^{-1}$ and~$\alpha''_i = \pi_\ell \circ \alpha_i \circ \iota_\ell^{-1}$. Then $((S_i\to S)_i, (\alpha'_i)_i)$ and $((S_i\to S)_i, (\alpha''_i)_i)$ are admissible pairs representing a level~$H$ structure $\alpha'$ and a level $B_0(\ell)$ structure~$\alpha''$ on $E$ respectively. Conversely, any pair $(\alpha',\alpha'')$ arises from a unique~$\alpha$ in this way.
\end{proposition}

\begin{proof}
    This is exactly the construction of \cite[§IV Construction 3.8]{Deligne.Rapoport1973}.
\end{proof}

\begin{theorem}\label{thm:algebraic diag}
    The moduli problems $\Pa_{H'}$ and $\Pa_H^{\ell\mathrm{-isog}}$ are isomorphic. In particular, the modular curve $Y_{H'}$ is a coarse moduli space for the problem $\Pa_H^{\ell\mathrm{-isog}}$. As a consequence, the diagram \cref{eq:diag pb} arises from morphisms of algebraic curves defined over~$\Q$
    \begin{center}\begin{tikzcd}
        & Y_{H'} \ar[ld] \ar[rd] & \\ Y_H & & Y_H,
    \end{tikzcd}\end{center}
    in the sense of \cref{thm:algebraic morphisms}.
\end{theorem}

\begin{proof}
By \cref{prop:moduli-product} and well-known results on the modular curve $Y_0(\ell)$ (see for instance \cite[Thm.~3.7.1, Cor.~3.7.2 and Cor.~6.8.7]{KatzMazur}), the moduli problems $\mathcal{P}_{H'}$ and $\mathcal{P}_H^{\ell\mathrm{-isog}}$ are isomorphic. Thus~\eqref{eq:moduli-diagram} yields natural transformations

\begin{displaymath}
    \begin{tikzcd}
        & \Pa_H' \ar[ld] 
        \ar[rd] & \\
        \Pa_H && \Pa_H,
    \end{tikzcd}
\end{displaymath}
Over~$\C$, these maps translate exactly into those from \cref{prop: maps over C}; they correspond to algebraic morphisms between the associated modular curves by \cref{rem:coarse-morphisms}.
\end{proof}

\end{document}